\documentclass[11pt]{article}
\usepackage{amsmath}
\usepackage{amsthm}
\usepackage{amssymb}
\usepackage{amscd}
\usepackage{color}
\usepackage[bookmarks=true]{hyperref}
\newcommand{\C}{\mathbb{C}}
\newcommand{\R}{\mathbb{R}}
\newcommand{\Z}{\mathbb{Z}}
\newcommand{\N}{\mathbb{N}}
\theoremstyle{plain}
\newtheorem{theorem}{Theorem}[section]
\newtheorem{proposition}[theorem]{Proposition}
\newtheorem{lemma}[theorem]{Lemma}
\newtheorem{Corollary}[theorem]{Corollary}

\theoremstyle{definition}
\newtheorem{definition}[theorem]{Definition}
\newtheorem{remark}[theorem]{Remark}

\makeatletter
\@addtoreset{equation}{section}
\makeatother
\begin{document}

%%%%%%%%%%%%%%%%%%%%%%%%%%%%%%%%%%%%%%%%%%%%%%%%%%%%%%%%%%%%%%%%%%%%%%%%%%%%%%%%%%%%%%%%%%%%%%%%%%%%%
%%%%%%%%%%%%%%%%%%%%%%%%%%%%%%%%%%%%%%%%%%%%%%%%%%%%%%%%%%%%%%%%%%%%%%%%%%%%%%%%%%%%%%%%%%%%%%%%%%%%%

\title{Some boundedness properties of solutions to the   
Vafa--Witten equations on closed four-manifolds 
} 

%%%%%%%%%%%%%%%%%%%%%%%%%%%%%%%%%%%%%%%%%%%%%%%%%%%%%%%%%%%%%%%%%%%%%%%%%%%%%%%%%%%%%%%%%%%%%%%%%%%%%

\author{  
Yuuji Tanaka}
\date{}

%%%%%%%%%%%%%%%%%%%%%%%%%%%%%%%%%%%%%%%%%%%%%%%%%%%%%%%%%%%%%%%%%%%%%%%%%%%%%%%%%%%%%%%%%%%%%%%%%%%%%%

\maketitle

%%%%%%%%%%%%%%%%%%%%%%%%%%%%%%%%%%%%%%%%%%%%%%%%%%%%%%%%%%%%%%%%%%%%%%%%%%%%%%%%%%%%%%%%%%%%%%%%%%%%%

\begin{abstract}
We consider a set of gauge-theoretic equations on closed oriented four-manifolds, 
which was introduced by Vafa and Witten. 
The equations involve a triple consisting of a connection and extra
 fields associated to a principal bundle over a closed oriented four-manifold.  
They are similar to Hitchin's equations over compact Riemann surfaces, 
and as part of the resemblance, there is no $L^2$-bound on the
 curvature without an $L^2$-bound on the extra fields.  
In this article, however, we observe 
that under the particular circumstance where the curvature 
does not become concentrated and the limiting connection is not locally reducible,
one obtains an $L^2$-bound on the extra fields. 
\end{abstract}

%%%%%%%%%%%%%%%%%%%%%%%%%%%%%%%%%%%%%%%%%%%%%%%%%%%%%%%%%%%%%%%%%%%%%%%%%%%%%%%%%%%%%%%%%%%%%%%%%%%%%%

%\renewcommand{\thefootnote}{\fnsymbol{footnote}}
%\footnote[0]{2000\textit{ Mathematics Subject Classification}.
% Primary 53C07 ; Secondary 49Q15, 53C15 ,  57R57 , 81T13.}

%%%%%%%%%%%%%%%%%%%%%%%%%%%%%%%%%%%%%%%%%%%%%%%%%%%%%%%%%%%%%%%%%%%%%%%%%%%%%%%%%%%%%%%%%%%%%%%%%%%%%%

%\markboth{}
%{}

%%%%%%%%%%%%%%%%%%%%%%%%%%%%%%%%%%%%%%%%%%%%%%%%%%%%%%%%%%%%%%%%%%%%%%%%%%%%%%%%%%%%%%%%%%%%%%%%%%%%%%

\section{Introduction}

In \cite{VW}, Vafa and Witten introduced a set of gauge-theoretic
equations on a 
four-manifold, motivated by the study of S-duality conjecture in a topologically twisted $N=4$ supersymmetric Yang--Mills theory. 
The moduli space of solutions to the equations is expected to produce
a possibly new invariant of some kind. 
Physicists believe that the partition functions which arise are the  
generating functions of the ``Euler characteristic'' of the ASD instanton moduli
spaces, 
and the S-duality conjecture predicts 
that the partition functions 
satisfy an enchanting modularity property
\cite{VW}. 
The equations were also considered by Haydys \cite{Hay} and Witten \cite{Wi} 
from a different point of view. 
They suggest a programme of 
``categorification'' of the Khovanov
homology by using a five-dimensional gauge theory, 
where the Vafa--Witten equations are a dimensional reduction of their
five-dimensional gauge-theoretic equation, and the solution spaces for
the Vafa--Witten equations are expected to be used in constructing a
Casson--Floer type theory.

The techniques developed by Taubes \cite{T1}, \cite{T2}, \cite{T3} and \cite{T5} 
have direct analogues in the 
Vafa--Witten case \cite{T4}, and we use them in this article.

The equations we consider involve a triple consisting of a
connection and other extra fields associated to a principal bundle over a 
four-manifold.  
Let $X$ be a closed, oriented, smooth Riemannian four-manifold with Riemannian
metric $g$, and let $P \to X$ be a principal $G$-bundle over $X$
with $G$ being a compact Lie group. 
We denote by $\mathcal{A}_{P}$ the set of all connections of $P$, 
and by $\Omega^{+} (X, \mathfrak{g}_{P})$ the set of self-dual two-forms
valued in the adjoint bundle $\mathfrak{g}_{P}$ of $P$. 
We consider the following equations for a triple $(A, B , \Gamma) \in
\mathcal{A}_{P} 
\times \Omega^{+} (X, \mathfrak{g}_{P}) \times \Omega^{0} (X,  \mathfrak{g}_{P})$, 
\begin{gather*}
d_{A} \Gamma + d_{A}^{*} B = 0, \\ 
F_{A}^{+}  + \frac{1}{8} [B.B ] + \frac{1}{2} [ B , \Gamma] =0 , 
\end{gather*}
where $[B. B ] \in \Omega^{+} (X, \mathfrak{g}_{P})$ is defined in
Section \ref{sec:vweq}. 
We call these equations the {\it Vafa--Witten equations}. 
We impose the additional simplifying assumption $\Gamma = 0$, 
which is not a significant restriction. 
As explained in Section \ref{sec:vweq} later, if $(A, B, \Gamma)$ is a solution, then so is $(A, B, 0)$. Moreover, the field $\Gamma$ vanishes if the connection is irreducible. 
Thus we consider only solution $(A,B)$ for which $\Gamma$ is assumed to vanish.

Mares studied analytic aspects of these equations in his
Ph.D thesis \cite{BM}. As in the case of Hitchin's equation \cite{Hi}, 
we do not have an $L^2$-bound on the curvature $F_{A}$ of a connection $A$
which satisfies the equations 
in general. 
In the case where the $L^2$-norm of $B$ is bounded, 
the following identity leads to a bound on the $L^2$-norm of $F_A$ for solutions.    
\begin{equation*}
\begin{split}
 \frac{1}{2} || d_{A}^{*} B ||_{L^2}^{2} 
   &+ || F_{A}^{+} + \frac{1}{8} [B . B] ||_{L^2}^{2} \\
  & = \int_{X} \left( \frac{1}{12} s |B|^2  
    - \frac{1}{2} W^{+} \langle B \odot B \rangle \right) \\ 
  & + \frac{1}{2} || F_{A} ||_{L^2}^{2} 
    + \frac{1}{4} || \nabla_{A} B ||_{L^2}^{2} 
    + \frac{1}{16} || [B . B] ||_{L^2}^{2}  
     + \frac{1}{2} \int_{X} \langle F_{A} \wedge F_{A} \rangle ,
\end{split}
\end{equation*}
where $\odot$ denotes some bilinear operation 
on $\Omega^{+} (X , \mathfrak{g}_{P}) \otimes \Omega^{+} (X ,
\mathfrak{g}_{P})$ (see \cite[\S B.4]{BM} for more details), 
$s$ is the scalar curvature of the metric, and  
$W^{+}$ is the self-dual part of the Weyl curvature of the metric. 
The last term on the right hand side of the identity is a topological
term.  
When $(A, B)$ is a solution, the left hand side vanishes. 
If $L^2$-norm of $B$ is bounded, 
then the integral in the middle line is bounded. 
Hence, using positivity of the middle two terms on the last line, 
one gets an upper bound on $|| F_{A} ||^{2}_{L^2}$.

Mares obtained the following compactness theorem, under the assumption
that the $L^2$ norm of $B$ is bounded. 
\begin{theorem}[\cite{BM}]
Let $\{ (A_{i} , B_{i}) \} $ be a sequence of solutions to the Vafa--Witten
 equations with $\int_{X} | B_{i} |^2  \leq C $. 
Then, by taking a subsequence of $\{ (A_{i} , B_{i}) \}$ if necessary, 
$\{ (A_{i} ,B_{i}) \}$ converges in $C^{\infty}$ on compact subsets
 to a smooth solution $(A , B)$ of the
 Vafa--Witten equations outside a finite set $\{ x_{1} , \dots , x_{k}
 \} \subset X$ after some gauge transformations on $X \setminus 
\{ x_{1} , \dots , x_{k}
 \}$. 
\label{th:BM}
\end{theorem}

In this article, we observe that there is a particular circumstance in which
one obtains an $L^2$-bound on the extra field $B$. 
The argument implies a non-constructive $L^2$-bound on $B \in \Omega^{+} ( X , \mathfrak{g}_{P})$
from the assumption that a limiting connection $A$ is not locally reducible. 
Here, we say that a connection $A$ on a principal $SU(2)$ or $SO(3)$ bundle $P$ is
{\it locally reducible} if the vector bundle $\mathfrak{g}_{P}$ has a 
$1$-dimensional subbundle that is $A$-covariantly constant (see
Definition \ref{def:locred} in Section \ref{sec:vweq}).

\begin{remark}
Note that a connection on a principal $SU(2)$ or $SO(3)$ bundle $P$ being locally reducible is the same as being honestly reducible if $X$ is simply-connected. 
\end{remark}

In order to properly state the main observation of this paper, we require
a technical assumption. 
We denote by $\delta$ the injectivity radius of $X$, so that geodesic
balls $B_r(x)$ of radius $r$ and centered at $x$ are well-defined
whenever $r<\delta$.

For a sequence of connections $\{ A_{i} \}$ on $P$, we define 
\begin{equation} 
S (\{ A_{i} \}) := \bigcap_{r \in (0, \delta)}\{ x \in X \, | \, \liminf_{i \to \infty} 
 \int_{B_{r} (x)} |F_{A_{i}}|^2  \, d \mu_{g} \geq \varepsilon_{\diamond}^2 \}, 
\label{eq:epdia}
\end{equation}
where $\varepsilon_{\diamond} >0$  is a positive constant (it is, in particular, independent of $\{ A_i \}$ and less than the constant of Uhlenbeck's gauge-fixing lemma). 
This number $\varepsilon_{\diamond}$ is defined as follows: 
Let $\kappa >1$ denote the number that appears in Proposition \ref{prop:tc}. 
Let $\varepsilon$ and $c$ denote the numbers appearing in Appendix \ref{app:bgb}'s equation \eqref{A1.7}. 
Set $\varepsilon_{\diamond} = \min(\kappa^{-1}, \varepsilon c^{-1}$). 
This set $S (\{ A_{i} \})$ describes the singular set of a sequence of
connections $\{ A_{i} \}$. 
In particular, if $x\not\in S$ then there is a small ball around $x$ for which some subsequence of connections has
curvature below the Uhlenbeck gauge fixing constant.
Our technical assumption is that 
$$ S (\{ A_{i} \}) = \emptyset .$$ 
Consequently, Uhlenbeck's theorem gurantees that there is some ball around every $x \in X$, and a subsequence such that the connections can be put into Coulomb gauge on this ball, with a uniform $L^2$-bound on the connection.

The main observation of this paper can be stated as follows.

\begin{theorem}
Let $X$ be a closed, oriented,  smooth Riemannian
 four-manifold 
with Riemannian metric $g$, 
and let $P \to X$ be a principal $G$-bundle over $X$ with $G$ being
 $SU(2)$ or $SO(3)$. 
Let $\{ ( A_{i} , B_{i} ) \}_{i \in \N}$ be a sequence of solutions to the
 Vafa--Witten equations with $S(\{ A_{i}\})$ being empty. 
Then, there exist a subsequence $\Xi \subset \N$ and a sequence of gauge
 transformations $\{ g_{i} \}_{i \in \Xi}$ such that $\{ g_{i} ( A_{i})
 \}_{i \in \Xi}$ converges weakly in the $L^{2}_{1}$-topology. 
If the limit is not locally reducible, 
then  there exists a positive number $C$ such that 
 $\int_{X} |B_{i}|^2 d \mu_{g} \leq C $ for all $i \in \Xi$, 
and $  \{ ( g_{i} ( A_{i}) , g_{i} (B_{i} ))
 \}_{i \in \Xi} $ converges in the 
 $C^{\infty}$-topology to a pair that obeys the Vafa--Witten equations. 
\label{th:main}
\end{theorem}

A modification of recent work of Taubes \cite{T2} and \cite{T3} is used
to prove this theorem.

As a particular case of Theorem \ref{th:main}, 
we have an $L^2$-bound on the extra fields in the ``fibre''  direction
at a connection
$A_{0}$ which is not locally reducible. 
Namely, 
\begin{Corollary}
Let $X$ be a closed, oriented, smooth Riemannian
 four-manifold with Riemannian metric $g$, 
and let $P \to X$ be a principal $G$-bundle over $X$ with $G$ being
 $SU(2)$ or $SO(3)$. 
Then, for any sequence of solutions  $\{ ( A_{0} , B_{i}) \}_{i \in \N}$ 
of the Vafa--Witten equations with $A_{0}$ being not locally reducible, 
there exist a subsequence $\Xi \subset \N$ and a positive constant $C >0$
 such that $\int_{X} |B_{i}|^2 d \mu_{g} \leq C $ for all $i \in \Xi$. 
\end{Corollary}

A proof of Theorem \ref{th:main} is given in Section 3. 
In Section 2, we describe some properties of the equations such as the
linearization, 
and also mention some remarks on reducible solutions to the equations.

In the rest of this section, let us mention an observation which led to
a formulation of Theorem \ref{th:main}. 
To make the argument rather simple, we assume that $X$ is
simply-connected so that a connection being locally reducible is in fact
honestly reducible. 
The point lies in relating the $L^2$-bound of a sequence of the extra
fields $\{ B_{i} \}$ to the irreducibility of the limiting connection.

Let $\{ (A_{i} , B_{i}) \}_{i \in \N}$ be a sequence of solutions to the
Vafa--Witten equations \eqref{VW1} and \eqref{VW2}, and suppose that $||
B_{i} ||_{L^2} \to \infty$ as $i \to \infty$. 
Put $r_{i} := || B_{i} ||_{L^2}$ and $\beta_{i} := B_{i} / r_{i} $ for
$i \in \N$. 
Then $(A_{i} , \beta_{i})$ satisfies 
$d_{A_{i}}^{*} \beta_{i} =0$, $F_{A_{i}}^{+}
 + \frac{1}{8} r_{i}^{2} [ \beta_{i} . \beta_{i}] =0$  
and $|| \beta_i ||_{L^2} =1$.

The technicalities of convergence of such a sequence 
are very delicate, and occupy the bulk of this paper.  But for now,
we assume that there is some hypothesis which ensures the convergence 
of $\{ (A_{i},\beta_{i}) \}$ to $\{ (A_{0},\beta)\}$ which satisfies 
$d_{A_{0}}^{*} \beta =0$,  
$[\beta . \beta ] =0$ and $|| \beta ||_{L^2} =1$.
After some further assumptions, this would imply that the connection 
$A_{0}$ is reducible by the following argument.

Firstly, we recall a notion of rank of a section $B \in \Omega^{+} (X, \mathfrak{g}_{P})$.

\begin{definition}
A $\mathfrak{g}_{P}$-valued self-dual two-form 
$B \in \Omega^{+} (X, \mathfrak{g}_{P})$ 
is said to be of {\it rank} $r$ if, when considered as a section of 
$\text{Hom} (\Lambda^{+, *} , \mathfrak{g}_{P})$, 
$B(x)$ has rank less than or equal to $r$ at every point $x \in X$ with
 equality at some point. 
\end{definition}

If the structure group of the principal bundle $P$ is either 
$SU(2)$ or $SO(3)$, 
then the possibilities for the rank of $B$ are $0, 1, 2$ or $3$. 
We next recall the following from \cite[\S~4.11]{BM}.

\begin{lemma}[\cite{BM}]
Let $P \to X $ be a principal $G$-bundle over a simply-connected
 Riemannian four-manifold $X$ with structure group $G$
 being either $SU(2)$ or $SO(3)$. 
Let $B \in \Omega^{+} (X, \mathfrak{g}_{P})$. 
If $[ B . B] =0$. then the rank of $B$ is at most one. 
\label{lem:BM}
\end{lemma}

Then, we invoke Lemma 4.3.25 in \cite{DK}.

\begin{lemma}[\cite{DK} Lem.~4.3.25]
Let $X$ be a simply-connected oriented Riemannian four-manifold, 
and let $P \to X $ be a principal $G$-bundle with structure group $G$
 being either $SU(2)$ or $SO(3)$. 
Let $\phi$ in $\Omega^{+} (X ,\mathfrak{g}_{P})$. 
Assume that $d_{A} \phi =0$. 
Then on any simply-connected open subset in X, where $\phi$ has rank 1 and is nowhere vanishing, the connection $A$ is reducible. 
\label{lem:DK}  
\end{lemma}

As $\beta \in \Omega^{+} (X, \mathfrak{g}_{P})$, it follows 
from $d_{A_0}^{*} \beta =0$ that $d_{A_0} \beta = 0$. 
From Lemma \ref{lem:BM}, $\beta$ is at most rank one. 
Thus, by Lemma \ref{lem:DK}, $A_0$ is reducible on an open subset in
$X$.  
We now further assume that $r_{i}^2 [ \beta_{i} . \beta_{i} ] \to 0$ so that $A_0$ is anti-self-dual. 
We then invoke the following.

\begin{lemma}[\cite{DK}, Lem.~4.3.21]
If $A$ is an irreducible $SU(2)$ of $SO(3)$ anti-self-dual connection on a bundle
 $E$ over a simply-connected four manifold $X$, 
then the restriction of $A$ to any non-empty open set in $X$ is also
 irreducible. 
\end{lemma}

From this, we deduce that $A_{0}$ is reducible on $X$ under the assumptions.

\paragraph{Acknowledgements.}
Some ideas of the proof of Theorem \ref{th:main} came out of discussions
with Dominic Joyce, while we confronted a similar phenomenon in the study
of Donaldson--Thomas instantons on compact symplectic 6-manifolds. 
I would like to thank him for sharing his insight and
expertise over the years. 
I would also like to thank Hirofumi Sasahira for valuable comments and
useful discussion on the subject.  
I am deeply grateful to  Cliff Taubes for pointing
out an error in the earlier version of this article, and providing me 
detailed notes, which clarify how to fix it by using part of his results on 
the Vafa--Witten analogue of \cite{T2} and \cite{T3}, 
and the so many corresponding enlightening discussions.  
I would also like to thank the referee for the many valuable 
comments and truly helpful suggestions.  
I am also grateful to Seoul National University, NCTS at National Taiwan University and BICMR at Peking University for support and hospitality, where part of this work was done during my visits in 2015--17. 
This work was partially supported by JSPS Grants-in-Aid for Scientific Research numbers JP15H02054 and JP16K05125.

\section{The Vafa--Witten equations on four-manifolds }
\label{sec:vweq}

In this section, we describe the precise form of the equations and  
their linearization, and then mention some remarks on reducible solutions.

\paragraph{Self-dual two-forms.}

For an oriented $n$-dimensional Riemannian manifold $X$ with Riemannian
metric $g$,  
there is an operator $* : \Lambda^{p} \to \Lambda^{n - p}$, 
called the {\it Hodge $*$-operator}, defined by 
$ \alpha \wedge * \beta = ( \alpha , \beta ) d \mu_{g}$, 
where $(,)$ denotes the natural metric on the forms defined by the
Riemannian metric, and $d \mu_{g}$ is the Riemannian volume element.

Setting $n=4$, the Hodge star operator on the space of two forms on $X$ induces the
following splitting into $\pm1$ eigenspaces: 
$$\Lambda^2 (X) = \Lambda^{+} (X) \oplus \Lambda^{-} (X). $$
Accordingly, the space of $\mathfrak{g}_{P}$-valued two-forms $\Omega^2 (X,
\mathfrak{g}_{P})$, where $P \to X$ is a principal $G$-bundle over $X$, 
splits as 
$$\Omega^{2} (X ,\mathfrak{g}_{P}) 
= \Omega^{+} (X ,\mathfrak{g}_{P}) \oplus \Omega^{-} (X ,
\mathfrak{g}_{P}). $$

The splitting of the space of two-forms can be viewed as that
corresponding to the isomorphism  
$\mathfrak{so}(4) \cong \mathfrak{so}(3) \oplus \mathfrak{so}(3)$ through the embedding 
$\Lambda^{2} T^{*}_{x}X \subset \text{End} (T^{*}_{x}X)$ as the
subalgebra of skew-symmetric endomorphisms. 
Hence $\Lambda^{+}$ is equipped with a pointwise
Lie-algebraic 
structure with the bracket denoted by $[\cdot , \cdot
]_{\Lambda^{+}}$. 
Using this together
with that of $\mathfrak{g}_{P}$, 
we define a bilinear map 
$[ \cdot .  \cdot ] : \left( \Lambda^{+} \otimes
  \mathfrak{g}_{P} \right)
\otimes  \left( \Lambda^{+}  \otimes \mathfrak{g}_{P} \right)  \to \Lambda^{+} \otimes
\mathfrak{g}_{P} $
by $\frac{1}{2} [ \cdot , \cdot ]_{\Lambda^{+}} \otimes [
\cdot , \cdot]_{\mathfrak{g}_{P}}$.

\paragraph{Equations.}

Let $X$ be a closed, oriented, smooth four-manifold, and let $P \to X$
be a principal $G$-bundle over $X$ with $G$ being a compact Lie group. 
We denote by $\mathcal{A}_{P}$ the set of all connections of $P$, 
by $\Omega^{+} (X, \mathfrak{g}_{P})$ the set of self-dual two-forms
valued in the adjoint bundle $\mathfrak{g}_{P}$ of $P$.

We consider the following equations for a triple $(A, B , \Gamma) \in
\mathcal{A}_{P} 
\times \Omega^{+} (X , \mathfrak{g}_{P}) \times \Omega^{0} (X , \mathfrak{g}_{P})$, 
\begin{gather}
d_{A} \Gamma + d_{A}^{*} B = 0, 
\label{VW01} \\ \
F_{A}^{+}  + \frac{1}{8} [B.B ] + \frac{1}{2} [ B , \Gamma] =0 . 
\label{VW02} 
\end{gather}
We call \eqref{VW01} and \eqref{VW02} the {\it Vafa--Witten equations}.

The left action of gauge group $\mathcal{G}_{P}$ is given by 
$( u(A) , \text{ad}(u)B , \text{ad} (u) \Gamma)$ 
for $u \in \mathcal{G}_{P}$ and  $(A, B , \Gamma) \in
\mathcal{A}_{P} 
\times \Omega^{+} (X, \mathfrak{g}_{P}) \times \Omega^{0} (X,
\mathfrak{g}_{P})$.

\paragraph{Linearization. }

Following \cite[\S 3.2.1]{BM}, let $(A,B,\Gamma)  \in
\mathcal{A}_{P} 
\times \Omega^{+} (X,  \mathfrak{g}_{P}) \times \Omega^{0} (X, 
\mathfrak{g}_{P})$ 
be a solution to the Vafa--Witten equations \eqref{VW01}
and \eqref{VW02}. 
Then, the infinitesimal deformation at $(A, B, \Gamma)$ is given by the following. 
\begin{equation}
\begin{split}
0 \longrightarrow  
\Omega^{0} ( X , \mathfrak{g}_{P}) 
\xrightarrow{d_{(A,B,\Gamma)}} 
&\Omega^{1} ( X , \mathfrak{g}_{P}) \oplus \Omega^{+} ( X , \mathfrak{g}_{P})
\oplus \Omega^{0} ( X , \mathfrak{g}_{P})  \\
 & \qquad \qquad 
 \xrightarrow{d_{(A,B,\Gamma)}^{+}} 
 \Omega^{+} ( X , \mathfrak{g}_{P}) \oplus \Omega^{1}
( X , \mathfrak{g}_{P}) 
\longrightarrow 0, \\
\label{eq:AHS}
\end{split}
\end{equation}
where the maps $ d_{(A,B,\Gamma)}$ and $ d_{(A,B,\Gamma)}^{+}$ are defined 
for $\xi \in \Omega^{0} ( X , \mathfrak{g}_{P}) $ 
and $(\alpha , \beta , \gamma) \in \Omega^{1} ( X , 
\mathfrak{g}_{P}) 
\oplus \Omega^{+} ( X , \mathfrak{g}_{P})
\oplus \Omega^{0} ( X , \mathfrak{g}_{P})$ by 
$$ d_{(A,B,\Gamma)} (\xi) 
:= ( - d_{A} \xi , [ \xi , B] , [ \xi , \Gamma] ) ,$$
\begin{equation*}
\begin{split}
& d_{(A,B,\Gamma)}^{+} 
( \alpha , \beta, \gamma )  \\
& \qquad := \left( d_{A}^{+} \alpha + \frac{1}{2} [B,\gamma] + \frac{1}{4} [B .\beta] - \frac{1}{2} [\Gamma
, \beta], d_{A} \Gamma - [\Gamma , \alpha] - [B , \alpha] + d_{A}^{*} \beta \right). \\  
\end{split}
\end{equation*}
The principal symbol of the operator $(d_{(A,B,\Gamma)}^{+} , 
d_{(A,B,\Gamma)}^{*})$ is self-adjoint, so the index of the complex
\eqref{eq:AHS} is always 
zero.

\paragraph{Local reducibility.}

\begin{definition}
A connection is {\it locally reducible} if there is an open cover of $X$ 
such that on each of the open subsets, there is a non-zero, covariantly constant section
 of $\mathfrak{g}_{P}$.
\label{def:locred}
\end{definition}

In the case considered here, the span of these local covariantly constant sections patch
together to define a covariantly constant line subbundle $\mathcal{I}$ in
$\mathfrak{g}_{P}$ on the
whole of $X$.  This implies in particular that the holonomy of the
connection lies in an $O(2)$ subgroup of $G$, which we are assuming $SO(3)$ or $SU(2)$.  Appendix \ref{app:redhol} to this paper reviews the relation between holonomy
and constant sections of associated bundles.
We also refer to the definition of irreducible connection explained in Appendix \ref{app:redhol}.

\paragraph{The equations in the irreducible semisimple case.}
As described in \cite[Th.2.1.1]{BM}, 
if $X$ is closed, the equations \eqref{VW01} and \eqref{VW02} are equivalent
to the following. 
\begin{gather*}
d_{A} \Gamma =0 , \, d_{A}^{*} B = 0,  \\ 
F_{A}^{+}  + \frac{1}{8} [B.B ] = 0 , \,    [ B , \Gamma] =0 . 
\end{gather*}
If $(A,B,\Gamma)$ is a solution, then it is clear from these equations that
$(A,B,0)$ is also a solution.  Furthermore,  
$\Gamma \in \Omega^0 (X , \mathfrak{g}_{P})$ lies in the kernel of 
$d_{A} : \Omega^{0} ( X , \mathfrak{g}_{P})
\to \Omega^1 (X , \mathfrak{g}_{P})$,
which is the infinitesimal stabilizer of $A$.   
If $G$ is semisimple, then the centre of $G$ is discrete, and so the
infinitesimal stabilizer of any irreducible $A$ is zero.   
In this case, $\Gamma$ vanishes. 
Since we are mainly interested in the irreducible semisimple case, we set $\Gamma=0$,
namely we consider only the following equations in the rest of this article. 
\begin{gather}
d_{A}^{*} B = 0 , 
\label{VW1} \\
F_{A}^{+} + \frac{1}{8} [B . B] =0 . 
\label{VW2} 
\end{gather}

\vspace{0.1cm}

We end this section by noting the following. If $A$ is a locally reducible connection on an
$SU(2)$ or $SO(3)$ bundle over $X$, then the adjoint bundle
$\mathfrak{g}_{P}$ 
decomposes as $\mathfrak{g}_{P} = \mathcal{I} \oplus C$, where
$\mathcal{I}$ is a real line bundle fixed by $A$, and $C$ is the
complementary $O(2)$ subbundle, as explained in Appendix \ref{app:redhol}. 
Hence, there are corresponding decompositions of 
$B = B_{\mathcal{I}} + B_{C}$ and $F_A^+=(F_A^+)_{\mathcal{I}}+(F_A^+)_C$, where 
the subscript $\mathcal{I}$ denotes projection onto 
$\Omega^+(X,\mathcal{I})$ and the subscript $C$ denotes
projection onto $\Omega^+(X, C)$.  
The connection $A$ is trivial on the $\mathcal{I}$ part,
so $(F_A^+)_{\mathcal{I}}=0$. 
According to \cite[\S 4.1.1]{BM}, this implies that $B_C$ has rank at most one, 
and hence $B$ has rank at most $2$.  
In this case, $[B.B]$  and hence $F_A^+$ have rank at most one.

\section{Proof of Main Theorem \ref{th:main}}

In this section, we show that Theorem \ref{th:main} is a corollary of the following theorem.

\begin{theorem}[\cite{T4}]
Let $\{ (A_{i} , B_{i})\}_{i \in \N}$ be a sequence of solutions to the
 Vafa--Witten equations with $S( \{ A_{i} \})$ being empty. 
Put $r_{i} := || B_{i} ||_{L^2} $ for $i \in \N$, 
and assume that
 $\{ r_{i} \}_{i \in \N}$ has no bounded subsequence. 
Then there exists a subsequence $\Xi \subset \N$ and a sequence of gauge
 transformation $\{g_{i} \}_{i \in \Xi}$ such that $\{ g_{i} (A_{i})
 \}_{i \in \Xi}$ converges in the weak $L^{2}_{1}$-topology on $X$ to a limit that is
 anti-self-dual and locally reducible.
\label{th:taubes} 
\end{theorem}

The author learned from Clifford Taubes \cite{T4} that one can prove Theorem \ref{th:taubes} 
by adjusting some of his results in \cite{T2} to the
situation here and also by using Theorem 1.2
of \cite{T3}.

\vspace{0.2cm}

\noindent 
{\it Proof of Theorem \ref{th:main}, assuming Theorem \ref{th:taubes}.}  
Suppose the hypotheses of Theorem \ref{th:main} are satisfied, namely that
the limit connection is not locally reducible.  We deduce firstly the $L^2$ bound on $B_i$,
and secondly the $C^\infty$ convergence.  The first proceeds by contradiction, 
assuming that the bound claimed in Theorem \ref{th:main} does not hold. Then the 
hypotheses of Theorem \ref{th:taubes} are satisfied, namely $\{r_i\}$ has no 
bounded subsequence.  The conclusion of
Theorem \ref{th:taubes} that the limit is locally reducible 
contradicts the hypothesis that the limit is not locally reducible in Theorem \ref{th:main}. 
Therefore the bound claimed in Theorem \ref{th:main} must hold.
For the second claim in Theorem \ref{th:main} of $C^\infty$ convergence, 
appeal to either Theorem \ref{th:BM} or a modification of the first bullet
of Proposition 2.1 in \cite{T2} while invoking that $S(\{ A_{i} \})$ is empty.
This implies the $C^{\infty}$-convergence of $\{ ( g_{i} (A_{i}) , g_{i} (B_{i})
)\}_{i \in \Xi}$ with a limit satisfying the Vafa--Witten equations
\eqref{VW1} and \eqref{VW2}.  
\qed

\vspace{0.2cm}

To complete the proof of Theorem \ref{th:main}, it remains only to prove Theorem \ref{th:taubes}.

\section[Proof of Theorem \ref{th:taubes}]{Proof of Theorem \ref{th:taubes}
\footnote{The details of the proof presented here were conveyed to the
author by Clifford Taubes.}}

In this section, we give the proof of Theorem \ref{th:taubes}. 
The proof must show that a limiting connection exists 
(this connection is denoted by $A_0$). 
This is the content of Proposition \ref{prop:uh}, whose proof is given in Appendix 
\ref{sec:proof3.4}.

\begin{proposition}[Weak $L^2_1$ compactness of non-concentrating connections modulo gauge \cite{T4}] 
Let $\{A_i \}_{i \in \N}$ denote a sequence of smooth connections on the principal $SU(2)$ bundle $P$ with $S(\{ A_i \}) = \emptyset$. 
Then, there exists a subsequence $\Pi \subset \N$ and a sequence of automorphisms $\{ g_i \}_{i \in \Pi}$ such that the sequence $\{ g_i (A_{i} ) \}_{i \in \Pi}$ converges weakly in the $L^2_1$ topology to a limit, which is a Sobolev class $L^2_1$ connection on $P$, 
denoted by $A_0$. 
\label{prop:uh}
\end{proposition}

By way of  a reminder, to say that $\{ g_i (A_{i} \}_{i \in \Pi}$ converges weakly in the $L^2_1$ topology to $A_0$ means the following: 
Fix a smooth connection on $P$ to be denoted by $\hat{A}$. 
Then write each $i \in \Pi$ version of $g_i (A_{i})$ as $\hat{A} + a_i$ with 
$a_{i}$ being a section of $T^* M \otimes \mathfrak{g}_P$ with $\mathfrak{g}_{P}$ denoting here (and subsequently) the Lie algebra of $SU(2)$. 
The weak $L^2_1$ convergence assertion requires that the elements in the sequence $\{ a_i \}_{i \in \Pi}$ have uniformly bounded $L^2_1$ norm and that this sequence converges weakly in the Hilbert space of $L^2_1$ sections of $T^*M \otimes \mathfrak{g}_{P}$. 
Let $a_0$ denote the limit $L^2_1$ section. 
Then $A_0 = \hat{A} + a_0$.

With $A_0$ in hand, it remains to prove that $A_0$ has anti-self-dual curvature tensor 
and that it is locally reducible. 
Local reducibility means, as per Definition \ref{def:locred}, that there is a rank one, 
real line subbundle $\mathcal{I} \subset \mathfrak{g}_{P}$, that is preserved by $A_0$-parallel transport. 
An equivalent formulation of local reducibility is as follows: 
There is a real line bundle $\mathcal{I}$ and a nonzero $A_0$-covariantly constant section (denoted by $\sigma_0$) of $\text{Hom} (\mathcal{I}, \mathfrak{g}_{P})$.  
The existence of $\mathcal{I}$ and $\sigma_0$ is asserted by Proposition \ref{prop:4.7}. 
The intermediate Propositions \ref{prop:taubes1}--\ref{prop:tc} and Theorem \ref{th:stZ} are needed for the construction of $\mathcal{I}$ and $\sigma_0$.

The task of proving local reducibility starts by considering the limit of the sequence 
$\{ \beta_{i} := B_{i} / r_i \}_{i \in \Pi}$ when $\{ r_i \}_{i \in \Pi}$ has no bounded subsequence. 
We show that away from the zero locus $Z$ of the limit, and up to gauge, a subsequence of $\beta_i$ converges to a tensor of rank one. 
The $\mathfrak{g}_{P}$ component of this tensor yields a section $\sigma_{\Delta}$ of $\text{Hom} (\mathcal{I} , \mathfrak{g}_{P})$ defined over $X \setminus Z$. 
Finally we show that $\sigma_{\Delta}$ extends over $Z$ to give the section $\sigma_0$ which is defined over all of $X$ and is $A_0$-covariantly constant, demonstrating that $A_0$ is locally reducible.

Firstly we have the following adjustment of Proposition
 2.1 in \cite{T2} to the Vafa--Witten case.

\begin{proposition}[\cite{T4}]
Assume that the sequence $\{ r_{n} \}_{n \in \N}$ has no bounded
 subsequence. 
Then there exists a subsequence $\Lambda \subset \N$ such that the
 following items a)--e) listed below hold. 
\begin{enumerate}
\item[a)]
The sequences $\{ \int_{X} ( \left| d |\beta_{i} | \right|^2 
 + |\beta_{i} |^2 ) \}_{i \in
	\Lambda}$ and $\{ \sup_{X} |\beta_{i} |\}_{i \in \Lambda}$ 
are bounded by a number that depends only on the
	geometric data. 
\item[b)] 
The sequence $\{ |\beta_{i} | \}_{i \in \Lambda}$ converges weakly 
in the $L^{2}_{1}$-topology and strongly in all $p < \infty$ version of
 the $L^{p}$-topology. 
The limit function is denoted suggestively by $| \hat{\beta}_{\diamond} |$,
this being a priori an $L^{\infty}$-function,
despite the fact that $\hat{\beta}_{\diamond}$ remains undefined as a 
section of any bundle.  
The function $| \hat{\beta}_{\diamond} |$ is defined at each point in $X$ by the
	rule whereby $ | \hat{\beta}_{\diamond} | (p) = \limsup_{i \in
	\Lambda} 
 | \beta_{i} | (p)$ for each $p \in X$. 
\item[c)] The sequence $\{ \langle \beta_{i} \otimes \beta_{i} \rangle
	\}_{i \in \Lambda}$ converges strongly in any $q < \infty$
	version of $L^q$-topology on the space of sections of
	$\Lambda^{+} \otimes \Lambda^{+}$. The limit section is denoted
	suggestively by $\langle \hat{\beta}_{\diamond} \otimes
	\hat{\beta}_{\diamond} \rangle $ and its trace is the function
	$| \hat{\beta}_{\diamond} |^2$. 
\item[d)] Use $f$ to denote any given $C^{0}$-function. 
\begin{enumerate}
\item[i)] The sequence $\{ r^{-2}_{i} \int_{X} f |F_{A_i}|^2 \}_{i \in \Lambda}$
      converges. 
\item[ii)] The sequences $\{ \int_{X} f | \nabla_{A_{i}} \beta_{i} |^2 \}_{i
      \in \Lambda}$ and $\{ 2 r_{i}^2 \int_{X} f | [\beta_{i} . \beta_{i}
      ] |^2 \}_{i \in \Lambda}$ converge. The limit of the first
      sequence is denoted by $Q_{\nabla , f}$ and that of the second by
      $Q_{\wedge , f}$. These are such that 
    $$ \frac{1}{2} \int_{X} (d^{*} d \, f) | \hat{\beta}_{\diamond} |^2  
      + Q_{\nabla , f} + Q_{\wedge , f} + \int_{X} f \, \mathcal{R} \left( 
   \langle \hat{\beta}_{\diamond} \otimes \hat{\beta}_{\diamond} \rangle
      \right) =0 , $$ 
    where $\mathcal{R}$ is an endomorphism associated to the Riemann curvature.  
\end{enumerate}
\item[e)] Fix $p \in X$ and let $G_{p}$ denote the Green's function with
	pole at $p$ for the operator $d^*d + 1$. The sequence
	$\{ \int_{X} G_{p} ( | \nabla_{A_{i}} \beta_{i} |^2 + 2 r_{i}^{2} |
	[\beta_{i} . \beta_{i} ] |^2 ) \}_{i \in \Lambda}$ is bounded by a number that
	depends only on the geometric data. 
     Let 
	$Q_{\diamond , p}$ denote the lim-inf of this sequence. The
	function $| \hat{\beta}_{\diamond} |^2$ obeys the equation 
         $$ \frac{1}{2} | \hat{\beta}_{\diamond} |^2 + Q_{\diamond , p }
	  = - \int_{X} G_{p} \left( \frac{1}{2} | \hat{\beta}_{\diamond}
	|^2 - \mathcal{R} ( \langle \hat{\beta}_{\diamond} \otimes
	\hat{\beta}_{\diamond} \rangle ) \right). $$
\end{enumerate}  
\label{prop:taubes1}   
\end{proposition}

The following is an analogue of Proposition 6.1 in \cite{T2}, 
but under the assumption that $S( \{ A_{i} \} )$ is empty, 
one does not need to follow the full argument of Section 5 in \cite{T2}.

\begin{proposition}[\cite{T4}]
The function $|\hat{\beta}_{\diamond}|$ of Proposition \ref{prop:taubes1} 
is continuous on $X$ and smooth at
 points where it is positive. 
Moreover, the sequence $\{ | \beta_{i} | \}_{i \in \Lambda}$ converges to
 $|\hat{\beta}_{\diamond}|$ in the $C^{0}$-topology on $X$.  
\end{proposition}

Section 3 of \cite{T2} introduces two positive numbers that measure 
the size of the curvature near any given point, these denoted by $r_{c
\wedge}$ and $r_{c F}$. 
Much of the analysis in Sections 3--6 in \cite{T2} is meant to deal with
the a priori fact that these numbers can be very small if the curvature
becomes concentrated near the chosen point. 
However, in the present context, the assumption that $S ( \{ A_i \}) = \emptyset$ in Theorems \ref{th:main}
and \ref{th:taubes} imply that both radii are bounded away from zero. 
Consequently, the assertions in Proposition 7.2 and Proposition 8.1 apply here. 
They are summarized as follows.

\begin{proposition}[\cite{T4}]
There exists $\kappa > 1$ with the following significance: 
Let $\{ (A_{i} , \beta_{i})
 \}_{i \in \N}$ be as described in Proposition \ref{prop:taubes1},  
let $\Lambda \subset \N$ be the subsequence from that proposition, 
and let $\delta$ denote the injectivity radius of $X$. 
Suppose that there exists $\rho \in ( 0 ,\delta )$ such that 
$$ \int_{B_p (\rho)}  
| F_{A} |^2 < \kappa^{-2} $$
for each $p \in X$ and $A \in \{ A_{i} \}_{i \in
 \Lambda}$. 
Then, there exist a closed set $Z \subset X$ with empty interior,  
a real line bundle $\mathcal{I}$ over $X \setminus Z$, 
a section $\nu \in \Gamma (X \setminus Z ; \mathcal{I} \otimes \Lambda^{+})$, 
a connection $A_{\Delta}$ on $P|_{X \setminus Z}$ 
and an isometric bundle homomorphism $\sigma_{\Delta} : \mathcal{I} 
\to \mathfrak{g}_{P}$. 
Their properties are listed below. 
\begin{itemize}
\item $Z$ is the zero locus of $|\hat{\beta}_{\diamond}|$. 
\item $| \hat{\beta}_{\diamond} |$ belongs to the H\"{o}lder space $C^{0, 1/ \kappa} (X)$.  
\item $|\nu|$ is the restriction of $|\hat{\beta}_{\diamond}|$ to $Z \subset X$. 
\item The section $\nu$ is harmonic in the sense that $d \nu =0$. 
\item $| \nabla \nu |$ is an $L^{2}$-function on $X \setminus Z$ that
      extends as an $L^{2}$-function on $X$. 
\item The curvature tensor of $A_{\Delta}$ is anti-self-dual. 
\item The homomorphism $\sigma_{\Delta}$ is $A_{\Delta}$-covariantly
      constant. 
There is, in addition, a subsequence $\Xi \subset \Lambda$ and a sequence
      $\{ g_{i} \}_{i \in \Xi}$ of automorphisms from $P$ such that 
\begin{enumerate}
\item[($i$)] $ \{ g_{i} ( \beta_{i}) \}_{i \in \Xi}$ converges to 
$\sigma_{\Delta} \circ \nu$ in the $L^{2}_{1}$-topology on compact subsets in
	     $X \setminus Z$ and in $C^{0}$-topology on $X$. 
\item[($ii$)] $\{ g_{i} ( A_{i}) \}_{i \in \Xi}$ converges on compact
	     subsets of $X \setminus Z$ to $A_{\Delta}$ in the
	     $L^{2}_{1}$-topology. 
\end{enumerate}
\end{itemize} 
\label{prop:tc}
\end{proposition}

The next theorem is a special case of Theorem 1.2 in \cite{T3}. 
It implies among other things that $Z$ has measure zero. 
(There are closed sets with empty interior and positive measure.)
This upcoming theorem uses the notion of a point of discontinuity
for $\mathcal{I}$; this being a point in $Z$ with the following
property: if $U$ is any neighbourhood of the point, then the restriction of $\mathcal{I}$ to $U \setminus Z$ is not isomorphic to the product bundle.

\begin{theorem}[\cite{T3}] 
Let $Z$ and $\mathcal{I}$ be as described in Proposition \ref{prop:tc}. 
The set $Z$ has Hausdorff dimension at most 2, and moreover, 
the set of the points of discontinuity for $\mathcal{I}$ (defined
in the preceding paragraph) are the points in
      the closure of an open subset of $Z$ that has the structure of a
      2-dimensional $C^1$-submanifold in $X$ denoted by $\Sigma$. 
\label{th:stZ}
\end{theorem}

We are finally ready to use the above results in the following proposition, of which items 1 and 3 conclude the proof of Theorem \ref{th:taubes}.

\begin{proposition}[\cite{T4}]
Let $\mathcal{I}$, $\sigma_{\Delta}$ and $A_{\Delta}$ be as described in Proposition \ref{prop:tc}, so that $\sigma_{\Delta}$ and $A_{\Delta}$ are defined over $X \setminus Z$. 
Then 
\begin{enumerate}
\item There exists a smooth anti-self-dual connection $A_0$ defined over all of $X$, and a continuous Sobolev-class $L_2^2$ gauge transformation $g_0$ defined over $X \setminus Z$, such that $g_0 (A_0)$ is the restriction to $X \setminus Z$ of $A_0$. 
Defining $\sigma_0 := g_0 ( \sigma_{\Delta})$ over $X \setminus Z$, then $\nabla_{A_0} \sigma_0 =0$. 
\item The bundle $\mathcal{I}$ over $X \setminus Z$ extends to a bundle defined over all of $X$, which we again denote by $\mathcal{I}$.  
\item There exist extensions of both $\nu \in \Gamma ( \mathcal{I} \otimes \Lambda^{+} )$ and $\sigma_0: \mathcal{I} \to \mathfrak{g}_{P}$ to all of $X$. 
We again denote these by $\nu$ and $\sigma_0$. 
The extensions satisfy $d \nu =0$ and $\nabla_{A_0} \sigma_{0} =0$.  
\end{enumerate}
\label{prop:4.7}
\end{proposition}

\begin{proof}

We first prove item 1.

Apply Proposition \ref{prop:uh} to the subsequence $\{ A_{i} \}_{i \in \Xi \subset \Lambda}$ of Proposition \ref{prop:tc}. 
This yields a subsequence $\Pi \subset \Xi$, a sequence of gauge transformations $\{ g_i \}_{i \in \Pi}$, and a connection $A_0$ which is the weak $L^2_1$ limit of
 $\{ g_i (A_{i} ) \}_{i \in \Pi}$ over all of $X$.

Recall from Proposition \ref{prop:tc} that $A_{\Delta}$ is the $L^2_1$ limit over compact subsets of $X \setminus Z$ of gauge transformations of $\{A_{i} \}_{i \in \Xi}$. 
In particular, both $A_0$ and $A_{\Delta}$ are weak $L^2_1$ limits over $X \setminus Z$ of gauge-equivalent connections. 
Since weak $L^2_1$ limits preserve $L^2_2$ gauge equivalence, it follows that there exists a Sobolev-class $L^2_2$ gauge transformation $g_0$ such that $g_0 (A_{\Delta}) = A_{0}$.

Note that $A_{\Delta}$ is anti-self-dual and gauge-equivalent over the complement of $Z$ to $A_0$. 
Thus, $A_0$ is an $L^2_1$ connection whose self-dual curvature is $L^2$, and vanishes on the complement of $Z$, which by Proposition \ref{th:stZ} is a set of measure zero. 
Hence the curvature of $A_0$ is anti-self-dual, and so a standard elliptic regularity argument (see e.g. \cite[\S 4.4]{DK} and Appendix \ref{sec:proof3.4}) can be used to prove that there is an $L^2_2$ and $C^0$ automorphism of $P$ that tranforms $A_0$ into a smooth connection with anti-self-dual curvature.  
After possibly composing $g_0$ with such an automorphism, we may assume without loss of generality that $A_0$ is smooth and that $g_0$ is continuous.

That $\nabla_{A_0} \sigma_0 = 0$ follows from Proposition \ref{prop:tc} since $\sigma_{\Delta}$ is $A_{\Delta}$-covariantly constant. This establishes the item 1.

We next prove the item 2, 
that $\mathcal{I}$ extends over $Z$. 
Let $\Sigma \subset Z$ denote the $C^1$ submanifold that is described by Theorem \ref{th:stZ}. 
It is enough to prove that $\Sigma$ is empty. 
For this purpose, assume to the contrary that $\Sigma \neq \emptyset$ and let $S \subset \Sigma$ be a component. 
This is a $C^1$ embedded, 2-dimensional surface. 
Fix a point $p \in S$. 
Since $S$ is $C^1$, there is an embedded disk $D \subset X$ whose closure intersects $S$ tranversally at a single point which is $p$. 
This is also its only intersection point with $Z$ since $S$ is an open subset of $Z$. 
Since $p$ is a point of discontinuity for the bundle $\mathcal{I}$, 
the restriction of $\mathcal{I}$ to $D \setminus p$ is not isomorphic to the product line bundle. 
In particular, parallel transport by $A_{\Delta}$ of $\sigma_{\Delta}$ along any circle in $D \setminus p$ which wraps once around $p$ gives $- \sigma_{\Delta}$. 
However, $A_{\Delta}$ is gauge-equivalent to a connection which is smooth over all of $D$. 
Thus parallel transport around sufficiently small circle in $D$ will be arbitrarily close to $+ \sigma_{\Delta}$, which is a contradiction.

Finally, we prove item 3 by showing that both $\nu$ and $\sigma_0$ extend to all of $X$ as $L^2_1$ sections. 
Granted this extension, we may argue as in item 1 that both 
$d \nu$ and $\nabla_{A_0} \sigma_0$ are $L^2$ sections which vanish almost everywhere, and hence by elliptic regularity, $\nu$ and $\sigma_0$ are smooth and satisfy $d \nu =0$ and $\nabla_{A_0} \sigma_0 = 0$ over all of $X$.

As a warning, note that not all $L^2_1$ functions extend over sets of measure zero. 
The simplest example is perhaps the Heaviside theta function, which is $L^2_1$ on $(-1, 0) \cup (0,1)$, but has no $L^2_1$ extension to $(-1, 1)$. 
Nevertheless, we will show that $Z$ is negligibly small in the sense that it has zero capacity and thus it follows that $L^2_1$ extensions over $Z$ always exist.

The section $\nu$ extends over $X$ as desired if there exists a
sequence $\{ u_{n} \}_{n \in \{ 1,2, \dots \}}$ of 
smooth functions on $X$ with pointwise norm bounded by 1 that
vanish on $Z$ and converge in the $L^2_1$ 
topology to the constant function 1. 
Indeed, given such a sequence, 
we may pass to a subsequence such that $\{u_n \}_{n \in \{ 1,2, \dots \}}$ converges to $1$ almost everywhere. 
By the dominated convergence theorem in $L^2$, it follows that $\{ u_n \nu \}_{n \in \{ 1,2, \dots \}}$ defines a Cauchy sequence of $L^2_1$ sections over the whole of $X$.  
The limit is equal to $\nu$ on $X \setminus Z$, and hence is our desired $L^2_1$ extension over $X$. 
The same argument also gives the desired $L^2_1$ extension of $\sigma_0$.

To find the desired sequence $\{ u_{n} \}_{n \in \{ 1,2, \dots \}}$ of
$L^2_1$ approximations to $1$, we first digress
to give a definition. 
Fix a ball $B \subset X$ of small radius and suppose that $V$ is a subset of $B$.
The capacity of the pair $(B,V)$ is denoted by $\text{cap}(B,V)$ 
and it is defined as follows:
$$ \text{cap} 
(B,V) = \inf_{u \in C^{\infty}(B)} 
\left\{ \int_{X} |du|^2 : u=1 \text{ on } X \setminus B \text{ and }
u=0 
 \text{ on } V \right\} . $$
For the present purposes, the significance of $\text{cap} (B,V)$ is that 
if it is zero, then there is a
sequence of functions on $X$ with $L^{\infty}$ norm $1$ that vanish on $V$ 
but converge to the constant function $1$ in $L^2_1 (X)$. 
Now, the definition is such that the capacity of a countable
union of sets is no greater than the sum of the capacities of the sets. 
See, e.g. Theorem 1.3 in \cite{R} by Reshetnyak or the
more accessible reference of Theorem 3.2 in \cite{C}
by Costea. Meanwhile, Lemma 3 in \cite{V} 
proves that a closed set in a 4-dimensional ball with finite 2-dimensional Hausdorff measure has zero
capacity. 
Granted these last two observations, it is enough to prove that $Z$ 
is contained in a countable union
of sets with finite 2-dimensional Hausdorff measure. 
\footnote{We know from Theorem \ref{th:stZ} only that the Hausdorff dimension of $Z$ is at most 2, but this does not guarantee that the Hausdorff measure is locally finite. 
For instance, the Hausdorff dimension of the graph of $\sin (1/x)$ has Hausdorff dimension 1.}
To see that this is the case, note that
because $\mathcal{I}$ extends over $Z$, 
it is locally isomorphic to the product bundle. As a
consequence, the section $\nu$ of $\mathcal{I} \otimes \Lambda^{+}$ 
extends to the whole of $X$ as a smooth, harmonic
section. One can then use Taylor's theorem in the manner of \cite{HHL} 
to see that $Z$ is
contained in a countable union of 2-dimensional, Lipshitz submanifolds. 
Since dimension
2, Lipshitz submanifolds have locally finite 2-dimensional Hausdorff
measure, 
it follows
that $\Sigma$ is contained in a countable union of 
submanifolds with finite 2-dimensional Hausdorff measure.
\end{proof}

\appendix

\section[Weak $L^2_1$ compactness of non-concentrating connections modulo gauge]{Weak $L^2_1$ compactness of non-concentrating connections modulo gauge \footnote{The proof presented here was also conveyed to the author by Clifford Taubes.}
}
\label{sec:proof3.4}

In this section we prove Proposition \ref{prop:uh}.

We begin with a brief outline of the proof idea.  
Cover $X$ by finitely many small balls $\{ B_r (p_i) \}_{i}$ of fixed radius $r$ such that Uhlenbeck's Theorem applies for a subsequence of $\{ A_{\alpha} \}$. 
For any fixed $\alpha$ in this subsequence, 
let $P_{\alpha}$ denote the principal bundle associated to a cocycle of trivializations over the $\{ B_{r} (p_i) \}_{i}$ which put $A_{\alpha}$ into Coulomb gauge. 
We assume that $A_{\alpha}$ is a smooth connection  so that $P_{\alpha}$ is a smooth principal bundle. 
Denote by $\eta_{\alpha}$ the canonical isomorphism of $P_{\alpha}$ with $P$. 
There is a subsequence such that 
in each ball, the pull-backs of the local connection forms converge weakly in $L^2_1$. 
There is a further subsequence such that the cocycles determining $P_{\alpha}$ converge weakly in $L^2_2$ to a $C^0$ cocycle which determines a bundle denoted by $P^{\diamond}$.

For $\alpha$ sufficiently large, so that the cocycle of $P^{\diamond}$ is sufficiently close to the cocycle of $P_{\alpha}$, there is an isomorphism $\phi_{\alpha} : P^{\diamond} \to P_{\alpha}$. 
Next we construct a family of smoothing $P^{\delta}$ of $P^{\diamond}$, parametrized by small positive $\delta$, and let $P'$ denote a particular choice of $P^{\delta}$. 
This comes with an isomorphism $\lambda : P' \to P^{\diamond}$. 
Finally, composing all these isomorphisms we get the composite isomorphism $$\Phi_{\alpha} : P' \xrightarrow{\lambda} P^{\diamond} \xrightarrow{\phi_{\alpha}} 
P_{\alpha} \xrightarrow{\eta_{\alpha}}  P. $$
These isomorphisms are constructed so that along the appropriate subsequence, 
$\Phi^{*}_{\alpha} (A_{\alpha})$ is weakly $L^2_1$ convergent in the smooth bundle $P'$. 
Choosing some fixed index $\alpha_0$, it follows that $\{ \Phi_{\alpha}\Phi^{-1}_{\alpha_0}\}$ is a sequence of automorphisms such that $\{ ( \Phi_{\alpha} \Phi^{-1}_{\alpha_0})^{*} (A_{\alpha}) \}$ is weakly $L^2_1$ convergent on the original bundle $P$, as desired.

\subsection{Background for analysis for the proof} 
\label{app:bg}

To set the notation for what is to come, suppose that $X$ is a
Riemannian 
4-manifold. Let $r_0 > 0$ be a small number chosen to be much less than the injectivity radius
of $X$. (We assume that $X$ is either compact, or if not, then the part of $X$ of interest has a
lower bound for its injectivity radius.) 
Gaussian coordinate charts can be defined on
balls of radius $r_0$ in $X$. Supposing that $p \in X$ and $r \in (0,
r_0)$, 
we use $B_r (p)$ to denote the
open, radius $r$ ball centered at $p$. If the point $p$ is immaterial to
the discussion, we will use
just $B_r$ to denote this ball. The boundary of the closure of $B_r$ is
denoted by $\partial B_r$. One
more piece of notation: we use $c_0$ to denote a number that is greater than 1 and independent
of any relevant choices made, such as a point $p$, radius $r$ or a given connection or other
geometric object. The precise value of $c_0$ can be assumed to increase between successive
appearances.
What follows next is a review of some of the basic analysis that is needed later.
Supposing that $p$ is a given point in $X$ and $r \in  (0, r_0)$, 
let $f$ denote a given Sobolev class $L^2_1$  function on $B_r$. 
This is to say that $f$ and its differential $df$ are square integrable on $B_r$.
The following are the basic function space inequalities that will be used:
\begin{itemize}
\item $\left( \int_{B_{r}} |f|^4 \right)^{\frac{1}{2}} \leq 
c_{0} \int_{B_r} \left( |df|^2 + r^{-2} |f|^2 \right) .$ 
\item $\int_{B_{r}} \frac{1}{\text{dist} (\cdot , p)^2} |f|^2 
\leq c_{0} \int_{B_r} \left( |df|^2 + r^{-2} |f|^2 \right) .$
\item 
If $f$ has compact support in $B_{r}$ or if the integral of $f$ over $B_r$ is zero, then the
preceding two inequalities hold without the $r^{-2} |f|^2$ term in the
integrand. 
\item $\int_{B_{r}} f^2 
\leq c_{0} \left( r^2 \int_{B_{r}} |df|^2 + r \int_{\partial B_{r}} f^2 \right).$
\vspace{-0.5cm}
\begin{equation}
\label{A1.1}
\end{equation}
\end{itemize}
\vspace{-0.4cm}
The first, third and fourth inequalities in this list are basic Sobolev inequalities. The
second inequality is called Hardy's inequality.

To set the stage for a second set of basic facts, suppose that $n$ is a positive integer
and let $B^{n} \subset \R^n$ denote the radius 1, open ball at the
origin. 
Now let $m$ denote a second positive
integer and $f$ denote a given smooth map from the closure of $B^n$ to $B^m$ mapping $0$ to $0$.
Suppose that $x$ is a Sobolev class $L^2_2$ map from $X$ to $B^n$.

\begin{flushleft}
{\sc Fact~1}: 
The composition of first $x$ and then $f$ defines a Sobolev class $L^2_2$
map from $X$ to $B^{m}$ whose $L^2_2$ norm is a priori bounded by $c_{f}$
 times 
the square of the $L^2_2$ norm of $x$ with $c_{f}$ being a number that depends only on $f$.
\end{flushleft}

\begin{flushleft}
{\sc Fact~2}: If $x$ is a Sobolev class $L^2_2$ map from $X$ to $B^{n}$ 
and if $\{ x_{\alpha} \}_{\alpha \in \{ 1,2 ,\dots \} }$  is a sequence of
Sobolev class $L^2_2$ maps from $X$ to $B^{n}$ with an a priori $L^2_2$ norm bound that converges
weakly in the $L^2_2$ topology to $x$, then the sequence 
$\{ f \circ x_{\alpha} \}_{\alpha \in \{1,2 ,\dots \}}$ converges weakly
 in the $L^2_2$ topology to $f \circ x$.
\end{flushleft}

The assertion made by {\sc Fact~1} can be proved by using the basic calculus to write
\begin{itemize}
\item $| ( f \circ x) | \leq c_f |x|$. 
\item $\nabla (f \circ x) = (df) |_{x} \cdot \nabla x . $
\item $\nabla (df \cdot dx ) = ( \nabla d f) |_{x} \cdot (dx \otimes dx ) + ( df )|_{x} \cdot
      \nabla dx .$
\vspace{-0.5cm}
\begin{equation}
\label{A1.2}
\end{equation}
\end{itemize}
\vspace{-0.4cm}

The desired bound on the $L^2$ norms of $f \circ x$ and its first derivatives follow directly from
the top two bullets in \eqref{A1.2}. The desired bound on the $L^2$ norm of the second derivatives
of $f \circ x$ follow from the third bullet of \eqref{A1.2} and the first bullet of \eqref{A1.1}.
To prove {\sc Fact~2}, first use {\sc Fact~1} to conclude that the
sequence $\{ f \circ x_{\alpha} \}_{\alpha \in \{1,2 , \dots \}}$ is
uniformly bounded in the $L^2_2$ topology for maps from $X$ to $B^{m}$. It follows as a
consequence that this sequence has a subsequence that converges weakly
in the $L^2_2$ topology. 
It is sufficient to show that $f \circ x$ is necessarily the weak limit and to this end, it
is sufficient to show that $f \circ x $ is the strong $L^2$ limit. 
That can be proved
by first using the fundamental theorem of calculus to write 
\begin{equation}
 (f \circ x ) - (f \circ x_{\alpha}) = ( x - x_{\alpha}) \int_{0}^{1} 
\nabla f |_{x + t (x -x_{\alpha})} dt 
\label{A1.3}
\end{equation}
and then invoking the fact that $\{ x_{\alpha} \}_{\alpha \in \{ 1,2 ,
\dots \}}$ converges strongly in the $L^2$ topology to $x$.

Uhlenbeck's paper \cite{Uh2} plays a central role. 
To review some of what is done in this
paper, fix a point $p \in X$ and a number $r \in (0, r_{0})$ 
The product $SU(2)$ bundle over $B_{r}$ is the
bundle $B_{r} \times SU(2)$. 
The corresponding product connection is denoted by $\theta_{0}$. 
Any connection on $B_{r} \times SU(2)$ can be written as $\theta_{0} +
\mathfrak{a}_{A}$ with $\mathfrak{a}_{A}$ being a section over $B_{r}$ of the
product bundle $B_{r} \times \mathfrak{su}(2)$. 
Here and below, we use $\mathfrak{su} (2)$ to denote the Lie algebra of
$SU(2)$. The curvature 2-form of $A$ is denoted by $F_{A}$. The group of automorphisms of the
bundle $B_{r} \times SU(2)$ is the space of maps from $B_{r}$ to
$SU(2)$. 
If $s$ is an automorphism of this
product bundle and $A$ is a connection, then $s$ can be used to pull back $A$ to get a new
connection. Writing $A$ as $\theta_{0} + \mathfrak{a}_{A}$ and $s^{*}A$ 
as $\theta_{0} + \mathfrak{a}_{s^{*}A}$, the $\mathfrak{su}(2)$-valued
1-forms 
$\mathfrak{a}_{A}$ and $\mathfrak{a}_{s^* A}$ are related by the rule
\begin{equation}
\mathfrak{a}_{s^{*}A} = s^{-1} \mathfrak{a}_{A} s + s^{-1} d s 
\label{A1.4}
\end{equation}
with $s$ here viewed as a map from $B_{r}$ to $SU(2)$.

The first basic result is Corollary 2.2 in \cite{Uh2}. A slightly weaker
version is stated next as a lemma.
\begin{lemma}
There exists a universal number (independent of $X$ and its metric) to be
denoted by $c$ which is greater than 1 and which has the following
 significance: 
Suppose that $p \in X$ and that $r \in (0, r_{0})$. 
Let $A$ denote a smooth connection on $B_{r}$ that obeys
$$ \int_{B_{r}} |F_{A}|^2 \leq c^{-2}$$
Then there exists an automorphism of $B_r \times SU(2)$ to be denoted by
 $s$ such that $s^{*}A$ 
can be written as $\theta_{0} + \mathfrak{a}_{s^{*}A}$ obeying
\begin{itemize}
\item $d^{*} \mathfrak{a}_{s^{*}A} =0 $. 
\item $x \cdot \mathfrak{a}_{s^{*}A} =0 $ on $\partial B_{r}$. 
\item $\int_{B_{r}} ( | \nabla \mathfrak{a}_{s^{*}A} |^2 
+ r^{-2} | \mathfrak{a}_{s^{*}A} |^2 ) \leq c_0^2 \int_{B_{r}} |F_{A}|^2 $.
\end{itemize} 
\label{lemA1.1}
\end{lemma}

Note that Uhlenbeck's Corollary 2.2 is stated using Sobolev class
      $L^2_{1}$ connections. We only need the version with the
      connection being smooth. 
(One can get Uhlenbeck's $L^2_1$ 
version in a straightforward way from the version above by taking
      limits.) 
By way of a
parenthetical remark, the number $c$ that appears in Lemma \ref{lemA1.1} can be
      chosen so that $\mathfrak{a}_{s^{*}A}$ 
also obeys a reverse of sorts to the third bullet's inequality:
\begin{equation}
\int_{B_{r}} |F_{A}|^2 \leq c_0^2 \int_{B_r} | \nabla \mathfrak{a}_{s^{*}A}
 |^2  . 
\label{A1.5}
\end{equation}
That this can be done follows from the integration by parts equality at the bottom of
page 35 in \cite{Uh2} and the inequalities in \eqref{A1.1}. 
The number $c$ is chosen in what follows so
that \eqref{A1.5} holds.
Let $p \in X$ and $r \in (0 , r_{0})$. 
Suppose now that $\{ A_{\alpha}\}_{\alpha \in \{ 1,2, \dots \}}$ is a sequence of
connections on $B_{r} \times SU(2)$, each with the $L^2$ norm on $B_{r}$ 
of the curvature of each $A_{\alpha}$ being
less than $c^{-1}$. Let $\{ s_{\alpha} \}_{\alpha \in \{ 1, 2, \dots
\}}$ 
denote the corresponding sequence of automorphisms of
$B_{r} \times SU(2)$ that is supplied by Lemma \ref{lemA1.1}. 
For each $\alpha \in \{ 1, 2, \dots \}$, write $s_{\alpha}^{*}
A_{\alpha}$ as $\theta_{0} + \mathfrak{a}_{\alpha}$ 
with $\mathfrak{a}_{\alpha}$ obeying the $A = A_{\alpha}$ 
version of the three bullets in Lemma \ref{lemA1.1}. 
In particular, it follows from the third bullet that the sequence 
$\{ \mathfrak{a}_{\alpha} \}_{\alpha \in \{ 1, 2, \dots \}}$ 
is uniformly bounded in the $L^2_1$ topology on $B_{r}$. 
Therefore, this sequence has a weakly convergent subsequence. 
Fix such a subsequence and let $\Theta$ denote the set of positive integers for the corresponding
indexing set in $\{ 1, 2, \dots \}$. 
The limit of $\{ \mathfrak{a}_{\alpha} \}_{\alpha \in \Theta}$ is a
Sobolev $L^2_1$ section over $B_{r}$ of $T^{*} B_{r} \times
\mathfrak{su}(2)$. 
Denote this limit by $\mathfrak{a}$. A priori, this limit 1-form obeys
\begin{itemize}
\item $d^{*} \mathfrak{a} =0 .$
\item $x \cdot \mathfrak{a} = 0$ on $\partial B_r$. 
\item $\int_{B_r} \left( | \nabla \mathfrak{a} |^2 + r^{-2}  |
      \mathfrak{a} |^2 \right) \leq 
 c_0^2 \liminf_{\alpha \in \Theta} \int_{B_r} | F_{A_{\alpha}}|^2 $. 
\vspace{-0.5cm}
\begin{equation}
\label{A1.6}
\end{equation}
\end{itemize}
\vspace{-0.4cm}
Note with regards to the second bullet that the components of the
restriction of $\mathfrak{a}$ 
to $\partial B_{r}$ 
are functions in the fractional Sobolev space $L^2_{1/2} (
\partial B_{r})$.

\subsection{Limits of sequences on $X$}
\label{app:bgb}

Let $P \to X$ denote a principal $SU(2)$ bundle. 
We assume the structure group of $P$ is just $SU(2)$ below, 
but a similar argument works for $SO(3)$ bundles as well. 
Suppose that $r \in ( 0 , r_0)$, that $N$ is a positive integer and that 
$\{ p_{i} \}_{ i \in \{ 1,2, \dots N \} }$ is a set of points in $X$ such that
the corresponding set $\{ B_{r} (p_{i}) \}_{i \in \{ 1,2, \dots N\}}$ 
is a cover of $X$ by balls of radius $r$.

Let $c$ be the constant defined in the previous section, 
and suppose that 
$\{ A_{\alpha }\}_{\alpha \in \{ 1,2, \dots \}}$ is a sequence of
connections on $P$ that obey the following condition: 
If $B_{r} \in \{ B_{r} (p_{i}) \}_{i \in \{ 1,2, \dots\}}$, 
then
\begin{equation}
\int_{B_{r}} |F_{A_{\alpha}} |^2 \leq  \varepsilon^2 c^{-2} .
\label{A1.7}
\end{equation}
The number $\varepsilon$ is determined near the end of this appendix.

Fix an index $i \in \{ 1, \dots N\}$ and let $B_{r}$ denote for the
moment the ball $B_{r} (p_i)$. Fix
an isomorphism $B_{r} \times SU(2)$ to $P|_{B_{r}}$. 
Denote this isomorphism by $\iota_{i}$. Invoke Lemma \ref{lemA1.1} 
for the sequence $\{ \iota_{i}^{*} A_{\alpha} \}$ on $B_{r}$. 
Denote the corresponding sequence of product
bundle automorphisms by $\{ s_{i, \alpha}\}_{\alpha \in \{ 1, 2, \dots\}}$ 
and for each index $\alpha$, write $s_{i ,\alpha}^{*} \iota_{i}^{*} A_{\alpha}$ as 
$\theta_{0} + \mathfrak{a}_{i ,\alpha}$. 
A subsequence in $\{ 1,2, \dots \}$ to be denoted by $\Theta$ can be chosen with the following
significance: For each index $i \in \{ 1, \dots , N\}$, 
the corresponding sequence $\{\mathfrak{a}_{i,\alpha}\}_{\alpha \in
\Theta}$ 
converges weakly in the $L^2_1$ topology. Use $\mathfrak{a}_{i}$ 
to denote the limit of this sequence. The
next lemma says that $\{ \mathfrak{a}_{i} \}_{i \in \{ 1, \dots , N \}}$
define a Sobolev class $L^2_{1}$ connection on a principal
$SU(2)$ bundle over $X$. By way of a reminder, a connection on a
principal bundle 
$P' \to X$ is said to be of Sobolev class $L^2_1$ 
if it can be written as $A_{0} + \hat{\mathrm{a}}$ with $A_{0}$ being a smooth
connection and with $\hat{\mathrm{a}}$ being an $L^2_1$ section of the vector bundle
$T^{*}X \otimes (P' \times_{SU(2)} \mathfrak{su}(2))$

\begin{proposition}
There is a smooth principal $SU(2)$ bundle $P' \to X$, a Sobolev class $L^2_1$
connection on $P'$ to be denoted by $A$, and, for each $i \in \{ 1,
 \dots N \}$, a $C^0$ and Sobolev class $L^2_2$ isomorphism from $B_{r}
 (p_{i}) \times SU(2)$ that pulls back $A$ as $\theta_{0} +
 \mathfrak{a}_{i}$. 
\label{propA1.2}
\end{proposition}

\begin{proof}
The proof has six steps.

\underline{Step 1}: Fix for the moment indices $i, k \in \{ 1, 2, \dots
 , N \}$ such that $B_{r} (p_{i}) \cap B_{r} (p_{k}) \neq \emptyset$.  
Use $g_{ki}$ to denote the automorphism $\iota_{k}^{-1} \circ \iota_{i}$
 of the product $SU(2)$ bundle over $B_{r} (p_{i}) \cap B_{r} (p_{k})$. 
For index $\alpha \in \Theta$, let $s_{ki, \alpha} = s_{k,\alpha }^{-1} 
 g_{ki} s_{i , \alpha}$. 
The definition is such that the two
connections satisfy 
$$\theta_{0} + \mathfrak{a}_{i,\alpha} 
= s_{ki,\alpha}^{*} (\theta_{0} + \mathfrak{a}_{k,\alpha}).  $$
Thus, 
\begin{equation*}
\begin{CD}
(B_r (p_i) \cap B_r (p_k) ) \times SU(2)  @>s_{i, \alpha}>>  (B_r (p_i) \cap B_r (p_k) ) \times SU(2) @>\iota_{i}>> P|_{B_r (p_i) \cap B_r (p_k)} \\
 @VVs_{ki, \alpha}V  @VVg_{ki}V  @VV\text{id}V \\ 
(B_r (p_i) \cap B_r (p_k) ) \times SU(2) @>s_{k, \alpha}>>  
 (B_r (p_i) \cap B_r (p_k) ) \times SU(2) @>\iota_{k}>> P|_{B_r (p_i) \cap B_r (p_k)} 
\end{CD}
\end{equation*}
commutes, and 
\begin{equation}
s_{ki,\alpha} \mathfrak{a}_{i ,\alpha} 
= 
\mathfrak{a}_{k ,\alpha} s_{ki , \alpha} 
+ ds_{ki, \alpha}.  
\label{A1.8}
\end{equation}
This equation implies among other things that the sequence $\{ s_{ki,
 \alpha} \}_{\alpha \in \Theta}$ is uniformly
bounded in the $L^2_2$ topology on the space of maps from $B_r(p_i) \cap
 B_r(p_k)$ to $SU(2)$. It
follows as a consequence that there is a subsequence of $\Theta$ such that the corresponding
subsequence of $\{ s_{ki, \alpha}\}_{\alpha \in \Theta}$ converges
 weakly in the $L^2_2$ topology to a Sobolev $L^2_2$ map from
$B_r(p_i) \cap  B_r (p_k)$ to $SU(2)$. We can assume without loss of generality that the same
subsequence works for all pairs $(i,k) \in \{ 1, \dots N \}$. We denote the latter sort of
subsequence by $\Lambda$. With this understood, let $s_{ki}$ denote the limit map.

\underline{Step~2}: As remarked above, the sequence $\{ s_{ki,
 \alpha}\}_{\alpha \in \Lambda} $ converges weakly in the $L^2_2$ topology
to an $L^2_2$ map $s_{ki}$ from $B_r (p_{i}) \cap B_r (p_k)$ to
 $SU(2)$. 
The next lemma lists three crucial
properties of the collection $\{ s_{ki} \}_{i,k \in \{ 1, \dots , N\}}$.

\begin{lemma}
The collection $\{ s_{ki} \}_{i,k \in \{ 1, \dots , N\}}$ has the
 following properties. 
\begin{itemize}
\item $s_{ki} \mathfrak{a}_{k} = \mathfrak{a}_{i} s_{ki} + d s_{ki}$ for
 each pair $(i,k)$ such that $B_{r} (p_i) \cap B_{r} (p_{k}) \neq
 \emptyset$. 
\item $ s_{ki} s_{ni}^{-1} s_{nk} = \mathbb{I}$ for each triple $(i ,k , n)$
 such that $B_{r} (p_{i}) \cap B_{r} (p_k) \cap B_{r(p_n)} \neq
 \emptyset$. 
\item $s_{ki}$ is {continuous} for each pair $(i,k)$ such that $B_{r}
 (p_i) \cap B_{r} (p_k) \neq \emptyset$. 
\vspace{-0.3cm}
\begin{equation}
\label{A1.9}
\end{equation}
\end{itemize}
\vspace{-0.4cm}
\label{lemA1.3}
\end{lemma}

\begin{proof}
The first bullet asserts an equality between $L^2_1$ functions and the
second asserts an equality between $L^2_2$ functions. They follow from the fact that the
sequences $\{ \mathfrak{a}_{i, \alpha} \}_{\alpha \in \Lambda}, 
\{ \mathfrak{a}_{k, \alpha} \}_{\alpha \in \Lambda}$ converge weakly in
 $L^2_1$ and the sequence $\{ s_{ki , \alpha } \}_{\alpha  \in \Lambda
 }$ converges weakly in $L^2_2$. 
Note in this regard that the product of a bounded, $L^2_2$ function
and an $L^2_1$ function is an $L^2_1$ function, and that the products of
 bounded $L^2_2$ 
functions are $L^2_2$ functions. 
However, if $s$ denotes the relevant
 bounded $L^2_2$ function, then the respective
norm bounds of these products depend on both the $L^2_2$ norm and $L^{\infty}$ 
norm of $s$.

The third bullet in \eqref{A1.9} is proved using the fact that $d^{*}
 \mathfrak{a}_{i} =0$ and $d^{*}
 \mathfrak{a}_{k} =0$. 
What with the first bullet in \eqref{A1.9}, these norm bounds imply that $s_{ki}$ obeys the equation
\begin{equation}
d^{\dagger} d s_{ki} = \langle \mathfrak{a}_{k} , d s_{ki} \rangle 
- \langle d s_{ki} , \mathfrak{a}_{i} \rangle ,
\label{A1.10}
\end{equation}
where $d^{\dagger}$ is short hand for $- * d *$  which is the formal
 $L^2$ adjoint of $d$; 
and where $\langle \alpha , \beta \rangle$ is
shorthand for $* ( \alpha \wedge * \beta)$.

To exploit this identity, fix an open set in $B_{r} (p_i) \cap B_{r}
 (p_{k})$ 
with compact closure in $B_{r} (p_i) \cap B_{r}
 (p_{k})$ to be denoted by $U_{ki}$. 
Then fix a smooth, non-negative function on
$B_{r} (p_i) \cap B_{r}
 (p_{k})$ to be denoted by $\chi_{ki}$ that is equal to 1 on a neighborhood of the set $U_{ik}$ and
0 on a neighborhood of the boundary of $B_{r} (p_i) \cap B_{r}
 (p_{k})$. Given a point $p \in B_{r} (p_i) \cap B_{r}
 (p_{k})$,
let $G_{p}$ denote the Dirichlet Green's function on $B_{r} (p_i) \cap B_{r}
 (p_{k})$  for the operator $d^{\dagger} d$ with
pole at $p$. Formally, $d^{\dagger} d G_{p} = \delta_p$ with
 $\delta_{p}$ 
being the Dirac delta-function with unit mass at $p$.
Note in any event that
\begin{equation}
|G_{p}| ( \cdot ) \leq c_{0} \frac{1}{\text{dist} (\cdot , p)^2} 
\, \text{ and } \, 
|\nabla G_{p}| \leq c_{0} \frac{1}{\text{dist} (\cdot , p)^3} 
\label{A1.11}
\end{equation}

It follows from \eqref{A1.10} that if $p \in U_{ki}$, then
\begin{equation}
\begin{split}
s_{ki} (p) 
&= \int_{B_{r} (p_i) \cap B_{r}
 (p_{k})} G_{p} \chi_{ki} 
( \langle \mathfrak{a}_{k} , d s_{ki} \rangle 
 - \langle d s_{ki} , \mathfrak{a}_{i} \rangle ) \\
& \qquad + \int_{B_{r} (p_i) \cap B_{r}
 (p_{k})}
G_{p} 
( (  d^{\dagger} d \chi_{ki} ) s_{ki}  - 
2 \langle d \chi_{ki} , d s_{ki} \rangle ) .
\end{split}
\label{A1.12}
\end{equation}
Note that the absolute value of the left most integral on the right hand side
 of \eqref{A1.12} 
can be bounded using the left most inequality in \eqref{A1.11}, Cauchy--Schwarz and the
 second bullet in 
\eqref{A1.1} by
\begin{equation}
\left| \int_{B_{r} (p_i) \cap B_{r}
 (p_{k})} G_{p} \chi_{ki} 
( \langle \mathfrak{a}_{k} , d s_{ki} \rangle - 
\langle d s_{ki} , \mathfrak{a}_{i} \rangle ) \right| 
\leq  c_{0} 
( || \mathfrak{a}_{i} ||_{*}
 +  || \mathfrak{a}_{k} ||_{*} )
|| d s_{ki} ||_{*} ,
\label{A1.13}
\end{equation}
where we introduce the notation whereby $|| \alpha  ||_{*}$ is the sum of
 the $L^2$ norm of $\nabla \alpha $ and the $L^4$ 
norm of $\alpha $ on $B_{r} (p_i) \cap B_{r} (p_k)$.  
Meanwhile, the right most integral in \eqref{A1.12} is no greater
than 
$$ c_{0} \, \text{dist} (p, \partial U_{ki})^{-3} 
( || d \chi ||_{L^1} + 
 || d^{\dagger} d \chi ||_{L^1} \, \text{dist} (p , \partial U_{ki} )) , $$  
this being a consequence of integration by parts, 
both inequalities in \eqref{A1.11} and the fact
that $s_{ki}$ has norm bounded by $c_0$.
The formula in \eqref{A1.13} is the key to proving that $s_{ki}$ is
 continuous. 
To say more, fix
for the moment $\rho \in ( 0 ,r)$ so that there is a ball of radius
 $\rho$ contained entirely in $U_{ki}$ with
the distance from any point in this ball to $\partial U_{ki}$ 
being greater than $100 \rho$. 
Let $B_{\rho}$ denote
such a ball. Suppose now that $p$ and $q$ are two points in $B_{\rho}$. 
Use \eqref{A1.12} to write
\begin{equation}
s_{ki} (p) - s_{ki} (q) = 
\mathcal{K}_{+} + \mathcal{K}_{-}
\label{A1.14}
\end{equation}
where the $\mathcal{K}_{+}$ is the contribution to the left most integral
 in \eqref{A1.12} 
from $B_{3 \rho}$ and where $\mathcal{K}_{-}$ is
the rest of this integral and the right most integral in \eqref{A1.12}. 
One finds using \eqref{A1.12} and
the second bullet of \eqref{A1.1} that
\begin{itemize}
\item $| \mathcal{K}_{+} | \leq c_{0} 
( || \mathfrak{a}_{i} ||_{* , 3 \rho} + 
|| \mathfrak{a}_{k} ||_{* , 3 \rho} ) || d s_{ki} ||_{* , 3 \rho} $. 
\item 
$|\mathcal{K}_{-} |
\leq || G_{p} - G_{q} ||_{L^{\infty}_{-} }
( || \chi ||_{L^2} ( || \mathfrak{a}_{i} ||_{*} + ||\mathfrak{a}_{k}||_{*} ) || d s_{ki} ||_{*} $ \\
\hspace{3cm} $ + || d^{\dagger} d \chi ||_{L^1} ||s_{ki}||_{L^{\infty} }
  + || d \chi ||_{L^{4/3}} || d s_{ki}||_{*} ) . $ 
\vspace{-0.5cm}
\begin{equation}
\label{A1.15}
\end{equation}
\end{itemize}
\vspace{-0.4cm}
Here, we define $|| \alpha ||_{* , 3 \rho}$ to denote the sum of the $L^2$
 norm of $\nabla \alpha $ over $B_{3\rho}$ and the $L^4$ norm 
of $\alpha$ 
over $B_{3 \rho}$ and $L^{\infty}_{-}$ to denote the $L^{\infty}$ norm over the complement of $B_{3 \rho}$.

With \eqref{A1.14} and \eqref{A1.15} understood, let $\delta > 0$ be arbitrary. 
Since $\mathfrak{a}_{i} , \mathfrak{a}_{k}$ and $d s_{ki}$ are all $L^2_1$
functions, there exists some $\rho_{\delta}$ such that their $L^2_1$ 
norms and $L^4$ norms on any ball of
radius $3 \rho_{\delta}$ or less is smaller than $\delta$. 
It follows as a consequence that if $p, q$ are in such a
ball, then $ | \mathcal{K}_{+} | \leq c_0 \delta^2$.  
Finally, once $\rho$ is fixed, we observe that $|| G_{p} - G_{q} ||_{L^{\infty}_{-} }
\leq c_0 \rho^{-2} \, \text{dist} (p,q)$, and so $| \mathcal{K}_{-}| \to 0$ as $\text{dist} (p, q) \to 0$. 
Since $|\mathcal{K}_{+}|$ and $|\mathcal{K}_{-}|$ can simultaneously be made arbitrarily small, this proves that $s_{ki}$ is continuous. 
\end{proof}

\underline{Step~3}: 
Since the set $\{ s_{ki}\}_{i,k \in \{ 1, \dots , N\}}$ 
are continuous maps from intersections of the
cover $\{ B_{r} (p_{i}) \}_{i \in \{ 1,2, \dots \}}$
to $SU(2)$ (the third bullet of Lemma \ref{lemA1.3}) and they obey the cocycle
condition (the second bullet of Lemma \ref{lemA1.3}), 
they a priori define a $C^{0}$ principal $G$-bundle
over $X$. Denote this bundle by $P^{\diamond}$ The bundle $P^{\diamond}$ 
has more than just a $C^0$ principal
bundle structure because its transition functions are also in the
 Sobolev space $L^2_2$. This
implies in particular that one can talk about $L^2_2$ sections of 
associated bundles to $P^{\diamond}$ and
$L^2_1$ connections. In particular, an $L^2_1$ connection 
is given by cocycle data consisting of a
collection of $L^2_1$ Sobolev class 1-forms with values in
 $\mathfrak{su}(2)$, 
one for each set in $\{ B_{r} (p_{i}) \}_{i \in \{ 1,2, \dots , N\}}$, 
with the collection subject to the cocycle condition that is depicted in the
top bullet of Lemma \ref{lemA1.3}.

\underline{Step~4}: Every $C^0$ principal bundle is isomorphic to some
 $C^{\infty}$ principal bundle, and
thus $P^{\diamond}$ is also. 
This being the case, the proof of Proposition \ref{propA1.2} is completed with an
argument to the effect that there is an isomorphism from some smooth principal bundle to
$P^{\diamond}$ that respects the $L^2_2$ structure; and in particular,
 that pulls back 
$L^2_1$ connections to $L^2_1$ connections. 
One way to construct the bundle and the desired isomorphism is to first
view $P^{\diamond}$ as the pull-back of the tautological $SU(2)$ bundle
 over 
$\mathbb{HP}^{N-1}$ with the latter being
the quaternionic projective space $(\times_{N} \mathbb{H} \setminus 0) /
 \mathbb{H}^{*}$, where $\mathbb{H}$ denotes the quaternions, $N$ is
 a large integer and $\mathbb{H}^{*}$
denotes the non-zero quaternions. 
To define the map, choose a partition of unity for the
cover $\{ B_r (p_{i}) \}_{i \in \{1 , \dots , N\}}$, 
this denoted by $\{ \varphi_{i}\}$. The map in question is denoted by
 $\Phi$. 
Given an index $i \in \{ 1, \dots , N\}$, and a point $p \in B_r
 (p_{i})$, 
the point $\Phi(p)$ is the equivalence class in
$\mathbb{HP}^{N-1}$ of a vector in $\mathbb{H}$ whose coordinates are as follows:
\begin{itemize}
\item The $i$'th coordinate of $\Phi(p)$ is $\varphi_{i} (p) 1$. 
\item If $k \in \{  1,  \dots , N \}$ and $k \neq i$ then the $k$'th
      coordinate of $\Phi(p)$ is 
\begin{enumerate}
\item[$a)$] $\varphi_{k} (p) s_{ki}(p) 1 $ \, if $ B_{r} (p_{i}) \cap B_{r}
	 (p_{k} ) \neq \emptyset $. 
\item[$b)$] $0$ \qquad \qquad \quad \, if $ B_{r} (p_{i}) \cap B_{r}
	 (p_{k} ) = \emptyset $. 
\end{enumerate}
\vspace{-0.5cm}
\begin{equation}
\label{A1.16}
\end{equation}
\end{itemize}
\vspace{-0.4cm}
This definition is consistent with regards to points in intersecting balls because of the
cocycle condition from the second bullet of Lemma \ref{lemA1.3}. 
The map $\Phi$ is observedly an $L^2_2$ 
and $C^0$ map because this is so for the function $\{ s_{ki} \}_{i,k \in \{
 1,2, \dots , N \}} $. The tautological $SU(2)$
bundle over $\mathbb{HP}^{N-1}$ is denoted in what follows by
 $\mathcal{P}_{N-1}$. 
It follows by construction that the
bundles $P^{\diamond}$ and $\Phi^{*} \mathcal{P}_{N-1}$ 
are given by the same cocycle data and so they are one and the
same bundle.

\underline{Step~5}: A standard way to smooth a $C^0$ map from one manifold to a second
compact manifold is to embed the second manifold in some finite dimensional
Euclidean space. One then smooths the embedding map as a map to the Euclidean space.
If the smoothing is sufficiently close to the original map, then its image will lie in a
tubular neighborhood of the embedding and so the tubular neighborhood projection will
compose with the smoothing to give a perturbation of the original map that is smooth.
Let $X$ denote the domain manifold. In particular, if the smoothing operation is defined by
convolution with a smoothing kernel with support very close to the
 diagonal in $X \times X$,
then the resulting map to the Euclidean space will be very close to the original.
Moreover, the smoothed map will also be very close to the original in a given Sobolev
topology if the original map is in the given Sobolev topology. (This
 follows from {\sc Fact~1}
in Appendix \ref{app:bg}.) 
Taking the image manifold to be $\mathbb{HP}^{N-1}$, one then obtains a smooth map
from $X$ to $\mathbb{HP}^{N-1}$ that is as close as desired to the
 original in the $C^0$ topology and also in
the $L^2_2$ topology. The pull-back of the tautological bundle by such a smoothing is
denoted here by $P^{\delta}$ with $\delta$ being an upper bound
 for the $C^{0}$ 
and $L^2_2$ distances between the
map $\Phi$ and the smoothing. 
Note in particular that $\delta$ can be as small as desired. Note also
that it follows from this construction that $P^{\delta}$ 
has a product structure on each ball $B_r (p_i)$ with smooth transition
 functions 
$\{ g^{\delta}_{ki} \}_{k,i \in \{ 1,2, \dots , N\}}$ 
that are $\delta$ close in the $C^0$ and $L^2_2$ topology
to those of $P^{\diamond}$. 
(This again follows from {\sc Fact~1} in Appendix \ref{app:bg}.)

\underline{Step~6}: Granted the preceding, one can invoke Proposition 3.2 in \cite{Uh2} and then
copy almost verbatim the proof of Corollary 3.3 in \cite{Uh2} to prove the following lemma:

\begin{lemma}
Fix a refinement $\{ U_{i}\}_{i \in \{ 1, \dots , N\}}$ of the cover 
$\{ B_{r} (p_i) \}$ such that each $U_{i}$ 
has compact closure in the corresponding $B_r (p_i)$. 
If $\delta$ is sufficiently small, then there is a
collection $\{ \rho_{i} : U_{i} \to SU(2) \}_{i \in \{ 1, 2, \dots \}}$ 
of $C^{0} \cap L^{2}_{2}$ maps such that $g^{\delta}_{ki} \rho_{i}
= \rho_{k} s_{ki}$ for all pairs of 
indices $i,k$ corresponding to intersecting balls. 
This collection $\{ \rho_{i} \}$ therefore defines a
$C^0$ and $L^2_2$ isomorphism from $P^\diamond$ to $P^{\delta}$.
\label{lemA1.4}
\end{lemma}

\begin{proof}
As noted above, the proof amounts to invoking Proposition 3.2 in
\cite{Uh2} and copying almost verbatim the arguments for Corollary 3.3
 in \cite{Uh2}. 
The only
changes are to take $p = 2$ and to use the a priori fact that the transition functions are close
in the $C^0$ topology. In Uhlenbeck's proof, she uses the fact that they are close in the
Sobolev topology to deduce $C^0$ closeness. Here, one has $C^0$
 closeness a priori. 
The basic
observation in this regard is that the exponential map on a Lie group is invertible on some
fixed neighborhood of the identity element. Moreover, the inverse on this neighborhood
has uniform bounds on its derivatives to any desired order. Thus, if one knows a priori
that a map has image in this neighborhood, then derivatives of the map can be bounded
knowing the derivative of its inverse under the exponential map, and vice-versa. (This
is, again, an application of {Fact~1} in Appendix \ref{app:bg}.)

We take $P'$ to be $P^{\delta}$ for any fixed $\delta$ which is sufficiently small to satisfy the hypothesis of Lemma \ref{lemA1.4}.  
Let $\lambda : P' \to P^{\delta}$ denote the corresponding isomorphism of Lemma \ref{lemA1.4}.  
We take $A$ to be the connection of $P'$ which corresponds under $\lambda$ to the connection on $P^{\diamond}$ determined by $\{ \mathfrak{a}_{i} \}_{i \in \{ 1, \dots, N \}}$. 
This concludes the proof of Proposition \ref{propA1.2}. 
\end{proof}
\end{proof}

\subsection{The convergence to the limit}
Proposition \ref{propA1.2} asserts that the limit data $\{
\mathfrak{a}_{i} \}_{i \in \{ 1, \dots , N \}} $ comes from an $L^2_1$ 
connection on a smooth principal bundle over $X$. To restate matters, there is a smooth
bundle $P' \to X$, an $L^2_1$ connection $A$ on $P'$, 
and, over each $B_r (p_i)$, an isomorphism $\tau_{i} : B_r(p_i) \times
SU(2) \to P'$ such that $\tau_i^* A = \theta_{0} + \mathfrak{a}_{i}$. 
We first check that $P'$ is isomorphic to the
bundle $P$ where the original sequence $\{ A_{\alpha} \}_{\alpha \in
\Lambda}$ lives. To see that this is so, we use the fact
that principal $SU(2)$ bundles are classified by their second Chern number which can be
computed using curvature integrals. In particular, the second Chern number of $P'$ is
\begin{equation}
- \frac{1}{4 \pi^2} \int_{X} \text{tr} (F_{A} \wedge F_{A}) ,
\label{A1.17}
\end{equation}
with the trace being the trace from the 2-dimensional representation of
$SU(2)$. 
Likewise, the second Chern class of $P$ is any $\alpha \in \Lambda$
version of 
\begin{equation}
- \frac{1}{4 \pi^2} \int_{X} \text{tr} ( F_{A_{\alpha}} \wedge
 F_{A_{\alpha}}) . 
\label{A1.18}
\end{equation}
To compute these numbers, we use the fact that $\text{tr} (F_{A} \wedge
F_{A})$ on a set $U_i$ from the refined
cover can be written as $ d( \text{tr}( \mathfrak{a}_i \wedge  d
\mathfrak{a}_i + \frac{1}{3} \mathfrak{a}_{i} \wedge \mathfrak{a}_{i} 
\wedge \mathfrak{a}_{i} ))$. 
A similar formula holds for the
integrand in \eqref{A1.18} on $U_i$, 
the only change being that $\mathfrak{a}_{i, \alpha}$, replaces
$\mathfrak{a}_{i}$. 
Now let $\{ \varphi_i\}_{i \in \{ 1, \dots , N\}}$ 
be a partition of unity subordinate to $\{ U_{i} \}_{i \in \{ 1,
\dots N \}}$. Then \eqref{A1.17} is equal to
\begin{equation}
- \frac{1}{4 \pi^2} \sum_{i \in \{ 1, \dots , N\}} 
\int_{X} \varphi_{i} d \text{tr} \left( \mathfrak{a}_{i} \wedge 
d \mathfrak{a}_{i}  + \frac{1}{3} \mathfrak{a}_{i} \wedge 
\mathfrak{a}_{i} \wedge \mathfrak{a}_{i} \right), 
\label{A1.19}
\end{equation}
which is, after an integration by parts,
\begin{equation}
\frac{1}{4 \pi^2} \sum_{i \in \{ 1, \dots , N\}} 
\int_{X} d \varphi_{i}  \text{tr} \left( \mathfrak{a}_{i} \wedge 
d \mathfrak{a}_{i}  + \frac{1}{3} \mathfrak{a}_{i} \wedge 
\mathfrak{a}_{i} \wedge \mathfrak{a}_{i} \right) .
\label{A1.20}
\end{equation}
A similar formula holds for \eqref{A1.18} with the 
replacement of $\mathfrak{a}_{i}$ by $\mathfrak{a}_{i , \alpha}$. 
Since $\varphi_{i}$ is smooth, 
and since $\{ \mathfrak{a}_{i,\alpha} \}_{\alpha \in \Lambda}$ 
converges weakly in $L^2_1$ and strongly in $L^2$, it follows that for each $i$,
the $\Lambda$-indexed sequences with $\alpha$ term given by
\begin{equation}
\int_{X} d \varphi_{i}  \text{tr} \left( \mathfrak{a}_{i, \alpha} \wedge 
d \mathfrak{a}_{i, \alpha}  + \frac{1}{3} \mathfrak{a}_{i,\alpha} \wedge 
\mathfrak{a}_{i,\alpha} \wedge  \mathfrak{a}_{i,\alpha} \right)
\label{A1.21}
\end{equation}
converges as $\alpha$ gets ever larger to the $i$'th term in
\eqref{A1.20}.

Now that we know that $P$ and $P'$ are isomorphic, we need to construct
a sequence $\{ g_{\alpha } \}_{\alpha \in \Lambda}$ 
of smooth isomorphisms from $P'$ to $P$ with the property that 
$\{ g_{\alpha}^{*} A_{\alpha} \}_{\alpha \in \Lambda}$ 
converges in the weak $L^2_1$ topology to $A$. 
The construction has seven steps.

\underline{Step~1}: 
To set the stage for this, fix a pair $i,k \in \{ 1, \dots , N \}$ 
whose corresponding balls intersect. Then fix an index $\alpha$. 
The transition function $s_{ki , \alpha}$ obeys \eqref{A1.8} and thus it
obeys the analogue of \eqref{A1.10}:
\begin{equation}
d^{\dagger} d s_{ki , \alpha } = 
\langle \mathfrak{a}_{k , \alpha} , d s_{ki , \alpha} \rangle 
- \langle d s_{ki, \alpha} , \mathfrak{a}_{i, \alpha} \rangle 
\label{A1.22}
\end{equation}

As a consequence, it is given on $U_{ki}$ by the analogue of \eqref{A1.12}:
\begin{equation}
\begin{split}
s_{ki, \alpha} (p) 
&= \int_{B_{r}(p_i) \cap B_{r} (p_{k})} 
 G_{p} \chi_{ki} ( \langle \mathfrak{a}_{k ,\alpha} , d s_{ki , \alpha}
 \rangle  - \langle d s_{ki, \alpha} , \mathfrak{a}_{i, \alpha} \rangle ) \\
&\qquad + \int_{B_{r}(p_i) \cap B_{r} (p_{k})} 
 G_{p} \left( - d^{\dagger} d \chi_{ki} s_{ki , \alpha } 
 + 2 \langle d \chi_{ki} , d s_{ki, \alpha} \rangle \right) . 
\end{split}
\label{A1.23}
\end{equation}
Now, if $p \in U_{ki}$, 
then the right most integral in \eqref{A1.23} converges uniformly as
$\alpha$ gets ever larger to the corresponding integral in \eqref{A1.12} 
because $G_{p}$ is smooth on the support of the
integrand as long as $p \in U_{ki}$. 
The left most integral in \eqref{A1.23} need not converge as $\alpha$ gets
larger. In any event, it follows from \eqref{A1.11} and from
\eqref{A1.1} 
using Lemma \ref{lemA1.1} that the absolute value of this integral is no greater
than $c_0 ( ||\mathfrak{a}_{k, \alpha} ||_{*} + || \mathfrak{a}_{i, \alpha} ||_{*} ) 
|| d s_{ki, \alpha} ||_{*}$. 
From \eqref{A1.8} 
and \eqref{A1.10} it follows that 
$$ ||d s_{ki, \alpha} ||_{*} \leq c_0 ( || \mathfrak{a}_{k ,\alpha} ||_{*} 
 + || \mathfrak{a}_{i, \alpha} ||_{*} )^2 . $$ 
Using \eqref{A1.7}, it follows that the left most integral in \eqref{A1.23} is bounded by $c_0 c^{-3} \varepsilon^{3}$.

\underline{Step~2}: 
Since this same bound holds for the norm of the left most integral in
\eqref{A1.12}, it follows that 
on the whole of $U_{ki}$, if $\alpha$ is sufficiently large, then 
\begin{equation}
| s_{ki, \alpha} - s_{ki} | \leq 
c_0 c^{-3} \varepsilon^3 . 
\label{A1.24}
\end{equation} 
This bound will determine the choice for $\varepsilon$.
If one can arrange for $\varepsilon$ to be small with a fixed size $N$ 
to the number of elements
in the cover $\{ B_{r} (p_{i} ) \}$, then one could invoke 
Uhlenbeck's Proposition 3.2 and the arguments of Corollary 3.3 
(using {Facts} 1 and 2 in Appendix \ref{app:bg}) to obtain the desired
sequence of isomorphisms. 
Unfortunately, our assumptions do not lead to an a priori
bound on the cover. Even so, one this difficulty can be surmounted to obtain the
desired conclusion. 
The next few steps explain how this is done.

\underline{Step~3}: 
The balls $\{ B_{r} (p_{i}) \}_{i \in \{ 1, \dots , N\}}$ 
can be chosen in any event so that at most $c_0$ of
them intersect at any given point, $c_{0}$ being independent of $N$. 
More to the point, one can
use Lemma 16 in \cite{D} 
 to partition $\{ 1 , \dots , N \}$ into
disjoint subsets $I_{1} \cup I_{2} \cup \cdots \cup I_{\ell}$ 
with $\ell \leq c_0$ such that if $i$ and $j$ are in the same subset, then
$B_r (p_i)$ has distance greater than $8r$ from $B_r(p_j)$. 
Construct such a partition, and then for
$\mathrm{a} \in \{ 1,2 , \dots , \ell\}$, 
let $\mathcal{U}_{\mathrm{a}}$ a denote the union of the balls $B_r (p_i)$ 
for $i \in I_{\mathrm{a}}$.

\underline{Step~4}: 
The collection $\{ \mathcal{U}_{\mathrm{a}}\}$ form an open cover of $X$ 
with the number of
elements being independent of $\varepsilon$ and the sequence $\{
A_{\alpha}\}$. 
The intersection between $\mathcal{U}_{\mathrm{a}}$ and
$\mathcal{U}_{\mathrm{b}}$ consist of the union of the intersections 
of the balls with index $i \in I_{\mathrm{a}}$ with those with
index $i \in I_{\mathrm{b}}$. 
Note that a given $B_{r} (p_{i})$ for $i \in I_{\mathrm{a}}$ 
intersects at most one ball $B_r (p_j)$ for $j \in I_{\mathrm{b}}$ 
because the distance between the balls labeled by $i \in I_{\mathrm{b}}$ 
is greater than the diameter of the balls. 
Thus, $\mathcal{U}_{\mathrm{a}} \cap \mathcal{U}_{\mathrm{b}}$ is the 
disjoint union of the sets $\{ B_{r} (p_i) \cap B_{r} (p_k) : 
i \in I_{\mathrm{a}} \text{ and } k \in I_{\mathrm{b} }\}$. 
There is one other point to note: If $\{ \chi_i \}_{i \in \{ 1, \dots ,
N\}}$ is a partition of unity subordinate to the
cover $\{ B_{r} (p_{i}) \}_{i \in \{ 1,2, \dots , N \}}$, 
then the collection
\begin{equation}
\left\{ \beta_{\mathrm{a}} = \sum_{i \in I_{\mathrm{a}}} \chi_{i}
\right\}_{\mathrm{a} \in \{ 1, \dots , \ell \}}
\label{A1.25}
\end{equation}
is a partition of unity subordinate to the cover
$\{ \mathcal{U}_{\mathrm{a}} \}_{ \mathrm{a} \in \{ 1, \dots , \ell \}}$.

\underline{Step~5}: Supposing that $\alpha \in \Lambda$, 
define a principal bundle over $X$ to be denoted by
$P_{\alpha}$ with $P_{\alpha}$ being defined by the following cocycle
data. It is given on any $\mathrm{a} \in \{ 1, \dots , \ell \}$ 
version of $\mathcal{U}_{\mathrm{a}}$ by $\mathcal{U}_{\mathrm{a}}
\times SU(2)$. 
If $i \in I_{\mathrm{a}}$ and $k \in I_{\mathrm{b}}$ 
are such that $B_{r} (p_{i}) \cap B_{r} (p_{k}) \neq \emptyset$, 
then $s_{ki, \alpha}$ 
is the transition on this component of $\mathcal{U}_{\mathrm{a}} \cap 
\mathcal{U}_{b}$. 
The bundle $P^{\diamond}$ also has a product structure on each set
$\mathcal{U}_{\mathrm{a}}$. If $i \in I_{\mathrm{a}}$ and 
$k \in I_{\mathrm{b}}$ are
such that $B_{r} (p_{i}) \cap B_{r} (p_{k}) \neq \emptyset$, then 
 $s_{ki}$ 
is the transition on this component of $\mathcal{U}_{\mathrm{a}} \cap 
\mathcal{U}_{b}$.

\underline{Step~6}: 
It follows from \eqref{A1.24} in the case when $\varepsilon$  is sufficiently small 
that Uhlenbeck's Proposition 3.2 and the arguments in her Corollary 3.3 
(with {Facts} 1 and 2 in Appendix \ref{app:bg})
can be invoked using the cover $\{ \mathcal{U}_{\mathrm{a}} 
\}_{\mathrm{a} \in \{ 1 , \dots , \ell\}}$ as input to find the following data when the
index $\alpha \in \Lambda $ is large: 
First, a refinement $\{ \mathcal{V}_{\mathrm{a}} \}_{\mathrm{a} \in \{
1 , \dots , \ell \}}$; and second, a set of continuous, $L^2_2$
maps $\{ \rho_{\mathrm{a} , \alpha} : \mathcal{V}_{\mathrm{a}} \to SU(2)
\}_{\mathrm{a} \in \{ 1, \dots , \ell\}}$ that obey 
\begin{itemize}
\item $s_{ki} \rho_{\mathrm{a} , \alpha} 
= \rho_{\mathrm{b} \alpha} s_{ki, \alpha}$ on $B_{r} (p_{i}) \cap B_{r}
      (p_{k})$ when $i \in I_{\mathrm{a}}$ and $k \in I_{\mathrm{b}}$
      and $B_{r} (p_{i}) \cap B_{r} (p_{k}) \neq \emptyset$. 
\item The $L^2_2$ and $C^{0}$ norms of $\{ \rho_{\mathrm{a},
      \alpha}\}_{\mathrm{a} \in \{ 1 , \dots , \ell \}}$ are less than
      $\varepsilon$ times an $\alpha$-independent number. 
\item The $L^2_1$ norms of $\{ \rho_{\mathrm{a} , \alpha} \}_{\mathrm{a}
      \in \{ 1 , \dots , \ell \}}$ have limit zero as $\alpha \to
      \infty$. 
\vspace{-0.5cm}
\begin{equation}
\label{A1.26}
\end{equation}
\end{itemize}
\vspace{-0.4cm}
Note in this regard that the third bullet follows from Uhlenbeck's
construction 
and {\sc Facts}
1 and 2 in Appendix \ref{app:bg} given that any $(i,k)$ version
of 
the sequence $\{s_{ki, \alpha}\}_{\alpha \in \Lambda}$ converges
strongly to $s_{ki}$ in the $L^2_1$ topology (because it converges weakly in
the $L^2_2$ 
topology.) In fact, there is strong convergence to zero in any
$L^{p}_1$ topology for $p < 4$.

\underline{Step~7}: Supposing that the index $\alpha \in \Lambda$ 
is large, then the collection $\{ \rho_{\mathrm{a} , \alpha
}\}_{\mathrm{a} \in \{ 1, \dots , \ell \}}$ 
define a $C^0$ and Sobolev class $L^2_2$ isomorphism 
$\phi_{\alpha} : P^{\diamond} \to P_{\alpha}$ 
with an a priori upper bound
on its $L^2_2$ norm. Meanwhile, 
the bundle $P_{\alpha}$ is, by definition, isomorphic to $P$ via an
isomorphism that identifies $A_{\alpha}$ on $B_r (p_i)$ for 
$i \in I_{\mathrm{a}}$ with $\theta_{0} + \mathfrak{a}_{i, \alpha}$. 
Let $\eta_{\alpha} : P_{\alpha} \to P $ denote the
latter isomorphism. 
Appendix \ref{app:bgb} defines a smooth principal bundle $P'$ with a $C^0$ and
Sobolev class $L^2_2$ isomorphism $\lambda : P' \to P^{\diamond}$ 
with the property that $\lambda^{-1}$ identifies Proposition
\ref{propA1.2}'s connection $A$ on $B_r (p_i)$ with $\theta_{0} +
\mathfrak{a}_{i}$. Let $\Phi_{\alpha}$ denote the isomorphism from $P'$ to $P$
that is obtained by composing first $\lambda$, then $\phi_{\alpha}$ 
and then $\eta_{\alpha}$. This is a $C^0$ and $L^2_2$ 
isomorphism with the property that $\{ \Phi_{\alpha}^{*} A_{\alpha}
\}_{\alpha \in \Lambda}$ is a sequence of $L^2_1$ connections that
converges weakly in the $L^2_1$ topology to the connection $A$. 
Since each $\Phi_{\alpha}$ is $C^0$ and $L^2_2$, we 
can use much the same argument as used in Steps 5 and 6 of the proof of
Proposition 
\ref{propA1.2} to smooth each large $\alpha$ version of $\Phi_{\alpha}$ 
slightly so that the result is a smooth isomorphism
from $P'$ to $P$ (to be denoted by $\Psi_{\alpha}$) 
such that the corresponding sequence $\Psi_{\alpha}^{*} A_{\alpha}$
converges weakly in the $L^{2}_1$ topology to the connection $A$. 
\qed

\section[Reducibility and holonomy]{Reducibility and holonomy 
\footnote{This appendix is due to the referee. The author
is grateful to the referee for allowing him  
to freely use it in this article. }}
\label{app:redhol}

In this appendix, we discuss the relation between reducibility of
connections on a principal $G$-bundle and holonomy in the cases $G=
SU(2)$ or $SO(3)$. 

As usual, we define the stabilizer $\Gamma_{A}$ of $A$ in the
gauge group $\mathcal{G}_{P}$ by 
$$ \Gamma_{A} := \{ u  \in \mathcal{G}_{P} \, | \, u (A) = A \} .$$
Given a connection $A$ on a principal $G$-bundle $P \to X$, there are various conflicting notions
of reducibility and irreducibility. Before giving our particular definitions for the cases $G =
SU(2)$ or $SO(3)$ in Definition \ref{def:red}, 
we review the standard definitions. According to \cite[\S 4.2.2]{DK}, a
connection $A$ is {\it reducible} if the holonomy of A is a proper
subgroup of $G$. 
It is {\it irreducible} if the stabilizer of $A$ under the group of gauge
transformations coincides with the center of $G$.

According to these definitions, when $X$ fails to be simply-connected, there are sometimes
examples of connections which are simultaneously reducible and irreducible. Such an example
is provided by any $SU(2)$ connection $A_0$ whose holonomy is a proper non-abelian subgroup
of $SU(2)$. 
Since the holonomy of $A_{0}$ is assumed to be a proper subgroup of
$SU(2)$, $A_0$ is certainly reducible. 
In order to show that it is also irreducible, we must show that the stabilizer of
$A_0$ coincides with the center $\mathbb{Z}_{2}$ of $SU(2)$. 
The stabilizer of $A_{0}$ is the centralizer of its holonomy. 
Being non-abelian, the holonomy is not contained in
any $U(1)$ subgroup. 
Thus the stabilizer of $A_{0}$ is a centralizer which is strictly smaller than the centralizer
of $U(1)$. The centralizer of any $U(1)$ subgroup of $SU(2)$ is the same
$U(1)$ subgroup. The only subgroups of $SU(2)$ which arise as
centralizers are $SU(2), U(1)$, and 
$\mathbb{Z}_{2}$. Since the only centralizer subgroup of $SU(2)$
smaller than $U(1)$ is the center $\mathbb{Z}_{2}$, it follows that $A_0$ must be irreducible.
Such a holonomy subgroup can not occur when $X$ is simply-connected, since in this case the holonomy
of any connection must be a closed and connected Lie subgroup. Thus when $G = SU(2)$ and X is simply-connected, there are no
intermediate holonomy subgroups between $U(1)$ and $SU(2)$.

The most useful definition of reducibility for Vafa--Witten theory over non-simply-connected
manifolds $X$ is the following, which we use in this article.
\begin{definition} 
A connection $A$ on a principal $SU(2)$ or $SO(3)$ bundle $P \to  X$ is
 {\it reducible} if
one of the following equivalent conditions is satisfied. 
\begin{itemize}
\item The stabilizer of $A$ under the group of gauge transformations has positive dimension.
\item  There exists a nonzero $\Gamma \in \Omega^0 (X ;
       \mathfrak{g}_{P})$ 
such that $d_{A} \Gamma =0$. 
\item The holonomy of $A$ is contained in some $SO(2)$ subgroup. 
\end{itemize}
A connection is said to be {\it irreducible} if it is not reducible. 
\label{def:red}
\end{definition}

Equivalence of the first two conditions in Definition \ref{def:red} 
follows from the fact that the Lie algebra of the stabilizer is the
kernel of $d_{A} : \Omega^{0} (X ; \mathfrak{g}_{P}) \to \Omega^{1} (X ;
\mathfrak{g}_{P})$. 
Equivalence of the last two conditions in Definition \ref{def:red}  
follows from the decomposition
$\mathfrak{g}_{P} = \underline{\R} \oplus C$, 
where $\underline{\R}$ is the span of $\Gamma$, and $C$ is the
complementary $SO(2)$ 
subbundle.

\begin{remark} 
When $G = SU(2)$, this definition of reducibility is equivalent to the holonomy of
$A$ being abelian. However, there are irreducible abelian connections when
$G = SO(3)$. 
For example, if $A$ preserves a splitting of $\mathfrak{g}_P$ as $\mathfrak{g}_{P} 
= \mathcal{I}_{1} \oplus \mathcal{I}_{2} \oplus (\mathcal{I}_{1} \otimes
 \mathcal{I}_{2})$ with each summand being a nontrivial real
line bundle, then the holonomy is $\mathbb{Z}_2 \times \mathbb{Z}_{2}$, 
which can not arise as a subgroup of $SO(2)$.
\end{remark}

With regards to the local reducibility introduced in Section
\ref{sec:vweq} of this article (Definition \ref{def:locred}), we have the following. 
\begin{proposition} 
A connection $A$ on a principal $SU(2)$ or $SO(3)$ bundle $P \to X$ is
 {\it locally
reducible} if and only if 
the holonomy of $A$ is contained in some $O(2)$ subgroup.
\end{proposition}

\begin{proof}
Every $O(2)$ subgroup of $SU(2)$ or $SO(3)$ is conjugate to the set of matrices of the form
\begin{equation*}
\left( 
\begin{matrix}
* & 0 \\
0 & * \\
\end{matrix}
\right)
\text{ and }
\left( 
\begin{matrix}
0 & * \\
* & 0 \\
\end{matrix}
\right)
\subset SU(2), 
\quad 
\left( 
\begin{matrix}
* & * & 0 \\
* & * & 0 \\
0 & 0 & * \\
\end{matrix}
\right)
\subset SO(3) , 
\end{equation*}
which is the subgroup such that the adjoint action preserves the subspace
\begin{equation*}
\sqrt{-1} \R
\left( 
\begin{matrix}
1 & 0 \\
0 & -1 \\
\end{matrix}
\right)
\subset \mathfrak{su}(2), 
\quad 
\R 
\left( 
\begin{matrix}
0 & -1 & 0 \\
1 & 0 & 0 \\
0 & 0 & 0 \\
\end{matrix}
\right)
\subset \mathfrak{so}(3). 
\end{equation*}
Whenever a 1-dimensional $A$-covariantly constant subbundle of
 $\mathfrak{g}_{P}$ exists, 
the holonomy must fix the subbundle, and hence must be contained in such
 a $O(2)$ 
subgroup. Conversely if holonomy is contained in an
$O(2)$ subgroup, parallel transport of the fixed subspace gives the
desired subbundle of $\mathfrak{g}_{P}$.
\end{proof}

\begin{remark}
If $A$ is locally reducible, then the restriction of $A$ to any simply-connected subset
of $X$ is reducible. This is because the holonomy subgroup over a simply-connected region
must be connected, and hence contained in an $SO(2)$ subgroup.
\end{remark}

\begin{remark}
If $A$ is a locally reducible $SU(2)$ or $SO(3)$ connection, then the following conditions
are equivalent.
\begin{itemize}
\item  $A$ is reducible.
\item One of the rank one subbundles $\mathcal{I}$ of $\mathfrak{g}_{P}$ 
which is preserved by $A$ admits a trivialization.
\item The first Stiefel--Whitney class vanishes for some rank one
      subbundle $\mathcal{I}$ which is preserved by $A$. 
\end{itemize}
Whenever $\pi_1 (X)$ has no subgroup of index two, it follows that
      $H^1(X; \mathbb{Z}_{2}) = 0$ 
and thus every locally reducible connection is reducible.
\end{remark}

%%%%%%%%%%%%%%%%%%%%%%%%%%%%%%%%%%%%%%%%%%%%%%%%%%%%%%%%%%%%%%%%%%%%%%%%%%%%%%%%%%%%%%%

\begin{flushleft}
Graduate School of Mathematics, Nagoya University, 
Furo-cho, Chikusa-ku, Nagoya, 464-8602, Japan \\
yu2tanaka@gmail.com
\end{flushleft}

%%%%%%%%%%%%%%%%%%%%%%%%%%%%%%%%%%%%%%%%%%%%%%%%%%%%%%%%%%%%%%%%%%%%%%%%%%%%%%%%%%%%%%%%%%%%%%%%%%%%%%

\end{document}